\documentclass{amsart} 
\usepackage{amsmath}
\usepackage{amssymb} 
\usepackage{amscd} 
\usepackage{graphicx}  
\usepackage{epsfig,comment}
\usepackage[all]{xy} 
\usepackage{marvosym}
\usepackage{float}
\usepackage{xcolor}
\usepackage{ascmac}
\usepackage[margin=1.2in]{geometry} 
\usepackage{hyperref}

\hypersetup{
	colorlinks=true,
	linktoc=all,
    	linkcolor={red!50!black},
   	citecolor={blue!50!black},
  	urlcolor={blue!80!black}
	}

\usepackage{booktabs} % Enhances quality of tables
\usepackage{tikz-cd,mathtools}

\usetikzlibrary{arrows}
\usetikzlibrary{braids}

\usepackage{wrapfig}
\usepackage{pdfpages}
\usepackage{mwe}

\usepackage{caption}

\DeclareMathOperator{\R}{\mathbb R}

\DeclareMathOperator{\Q}{\mathbb Q}
\DeclareMathOperator{\Z}{\mathbb Z}

\DeclareMathOperator{\brl}{brloc}
\DeclareMathOperator{\wrap}{wrap}

\DeclareMathOperator{\Hom}{Hom}

\DeclareMathOperator{\dc}{dc}

\newcommand{\mc}{\mathcal}

\theoremstyle{plain}
\newtheorem{theorem}{Theorem}
\newtheorem{lemma}[theorem]{Lemma}
\newtheorem{proposition}[theorem]{Proposition}
\newtheorem{corollary}[theorem]{Corollary}
\newtheorem{conjecture}[theorem]{Conjecture}

\newtheorem{question}[theorem]{Question}
\newtheorem{fact}[theorem]{Fact}
\theoremstyle{definition}
\newtheorem{remark}[theorem]{Remark}

\newtheorem{definition}[theorem]{Definition}
\newtheorem{notazione}[theorem]{Notation}
\newtheorem{example}[theorem]{Example}
\newtheorem{exercise}[]{Exercise}
\newtheorem{problem}[theorem]{Problem}

\newtheorem{vuoto}[theorem]{}

\numberwithin{theorem}{section}

\newcommand{\bt}{\begin{theorem}}
\newcommand{\et}{\end{theorem}}

\newcommand{\bv}{\begin{vuoto}}
\newcommand{\ev}{\end{vuoto}}

\newcommand{\bl}{\begin{lemma}}
\newcommand{\el}{\end{lemma}}

\newcommand{\bd}{\begin{definition}}
\newcommand{\ed}{\end{definition}}

\newcommand{\beq}{\begin{equation}}
\newcommand{\eeq}{\end{equation}}

\newcommand{\bexa}{\begin{example}}
\newcommand{\eexa}{\end{example}}

\newcommand{\bexe}{\begin{exercise}}
\newcommand{\eexe}{\end{exercise}}

\newcommand{\bfact}{\begin{fact}}
\newcommand{\efact}{\end{fact}}

\newcommand{\bprop}{\begin{proposition}}
\newcommand{\eprop}{\end{proposition}}

\newcommand{\bp}{\begin{proof}}
\newcommand{\ep}{\end{proof}}

\newcommand{\bc}{\begin{corollary}}
\newcommand{\ec}{\end{corollary}}

\newcommand{\bq}{\begin{question}}
\newcommand{\eq}{\end{question}}

\newcommand{\bcong}{\begin{conjecture}}
\newcommand{\econg}{\end{conjecture}}

\newcommand{\bproblem}{\begin{problem}}
\newcommand{\eproblem}{\end{problem}}

\newcommand{\bs}{\begin{proof}[Proof.]}
\newcommand{\es}{\end{proof}}

\newcommand{\br}{\begin{remark}}
\newcommand{\er}{\end{remark}}

\newcommand{\bn}{\begin{notazione}}
\newcommand{\en}{\end{notazione}}

\begin{document}
\baselineskip 13pt

\title{Sutured manifold hierarchies and the Thurston norm}
\author{Alessandro V. Cigna}
\address{Department of Mathematics, King's College London, Strand, London, UK, WC2R 2LS}
\email{alessandro.cigna@kcl.ac.uk}


\begin{abstract}
    \noindent By the classical work of Thurston and Gabai, the Thurston norm of any compact oriented irreducible $3$-manifold with toroidal boundary is determined by finitely many taut sutured manifold hierarchies. We give an explicit procedure to extract this information from such hierarchies. This is achieved via the maw dual graph construction, which can be incorporated into a general method for computing the Thurston norm of a manifold. As an application, we compute the Thurston norm of the exterior of all alternating and some non-alternating pretzel links with three components. Using these computations, we give a negative answer to a question of Baker--Taylor. Namely, we exhibit an infinite family of manifolds, and an infinite family of knots in each manifold, such that the wrapping number associated with these knots is not a seminorm on the respective second real homology group.

   Next, we show that if a taut surface in a link exterior is disjoint from a boundary torus and does not represent a corner of the Thurston norm ball, then the corresponding link component is algebraically split from the rest of the link. Moreover, the surface remains norm-minimising under every non-longitudinal Dehn filling on that boundary component and, under an assumption of atoroidality, is a leaf of a taut foliation transverse to the filling core.
\end{abstract}

\maketitle

\section{Introduction}
Given a Seifert surface $S$ for a link $L$, when does $S$ realise the minimal genus among Seifert surfaces for $L$? More generally, given a properly embedded oriented surface $S$ in a compact oriented $3$-manifold $M$, when does $S$ attain the minimal topological complexity among surfaces homologous to $S$? Here, our measure of ``topological complexity" is the Thurston norm (see Section \ref{sec: thurston norm} for details). In 1983, Gabai gave a complete answer to these questions in terms of sutured manifold hierarchies. 

\bt\cite{gabai_foliations_1983}\label{thm: gabai intro} Let $M$ be a compact oriented and irreducible $3$-manifold $M$ with (possibly empty) toroidal boundary and let $S\subset M$ be a properly embedded oriented surface. Suppose that $S$ meets each component of $\partial M$ in parallel oriented curves and no union of components of $S$ is zero in homology. Then, $(M,\partial M)$ admits a taut sutured manifold hierarchy starting with a decomposition along $S$ if and only if $S$ is norm-minimising for its class in $H_2(M,\partial M)$. 
\et

Gabai's theorem guarantees that taut sutured manifold hierarchies determine the Thurston norm of any integral class $\alpha\in H_2(M,\partial M)$. In Gabai's theory, taut foliations play the connecting bridge between the Thurston norm and sutured manifold hierarchies. Indeed, Gabai showed that taut sutured manifold hierarchies carry taut foliations and, for such foliations, Thurston's Inequality holds.

\bt[Thurston's Inequality \cite{Thurston1986ANF}] Let $\mc F$ be a taut foliation of a compact oriented $3$-manifold $M$, with $\mc F$ transverse to $\partial M$. For every properly embedded oriented surface $S\subset M$ without sphere nor disk components, the following inequality holds:
\begin{equation*}
    |\langle e(\mc F), [S]\rangle| \le |\chi(S)|.
\end{equation*}
Here, $e(\mc F)\in H^2(M,\partial M)$ is the relative Euler class of $\mc F$ and $[S]\in H_2(M,\partial M)$. Moreover, the equality is achieved if $S$ is a union of consistently oriented leaves of $\mc F$.
\et

In particular, each taut foliation of $M$ induces a lower bound on the topological complexity of \emph{any} surface $S \subset M$ without sphere or disk components, and this bound depends only on its homology class. In the present article, we highlight that a taut sutured manifold hierarchy indeed provides \emph{global} information about the Thurston norm and describe how to obtain this information. The method is theoretically simple: if we want to show that a surface $S$ is norm-minimising, we have to find a taut sutured manifold hierarchy for $(M,\partial M)$ starting with a decomposition along $S$. Once such a hierarchy has been found, we can use it to compute the relative Euler class of any taut foliation $\mc F$ carried by it. At this point, Thurston's Inequality restrains the Thurston unit ball to lie between the two hyperplanes $\{\alpha\in H_2(M,\partial M;\R)\, |\; |\langle e(\mc F),\alpha\rangle|=\pm 1\}$. Moreover, if surfaces $S_1,...,S_n$ are algebraically fully marked with respect to $\mc F$, meaning that $\langle e(\mc F),[S_i]\rangle=\chi(S_i)$ for each $i$, then also $S_1,...,S_n$ are norm-minimising and their homology classes project to the same face of the Thurston unit ball (see Corollaries \ref{cor: thurston} and \ref{cor: norm detection criterion}). 

The above strategy can be incorporated into a general approach to compute the Thurston norm of a compact, orientable $3$-manifold. We can synthetically outline its steps as:
 \vspace{-1pt}
\begin{align*}
&\hspace{-4pt}\boxed{\mathrm{candidate\ surfaces\ } S,S_1, ..., S_n}\\ &\hspace{56pt}\downarrow \\  &\hspace{28pt}\boxed{\mathrm{hierarchy\ for\ } S}\\&\hspace{58pt} \bigg\downarrow [\mathrm{Gab87b}] \\  &\hspace{25pt}\boxed{\mathrm{branched\ surface}} \\ &\hspace{58pt}\bigg\downarrow [\mathrm{Cig26}] \\ &\hspace{36.5pt}\boxed{\mathrm{Euler\ class}} \\ &\hspace{58pt}\bigg\downarrow[\mathrm{Thu86}] \\ &\hspace{27pt}\boxed{\mathrm{Thurston\ norm}}
\end{align*}

By \cite[Construction 4.16]{gabai_foliations_1987}, the whole information on a hierarchy can be reorganised in a branched surface carrying the foliation. The step of computing the Euler class of a foliation, given a branched surface $\mc B$ carrying it, is addressed in \cite{maw_dual_graph}. The idea is to use the dual graph of the branched surface to define a simplicial $1$-cycle called the \emph{maw dual graph} $\Gamma_m(\mc B)$. The graph $\Gamma_m(\mc B)$ represents the Poincar\'e dual of the relative Euler class of any foliation carried by $\mc B$. We emphasise that the only real novelty in the proposed method is the passage from hierarchies to Euler classes via the maw dual graph.
In Section \ref{sec: pretzels}, we illustrate this method while computing the Thurston norm in the exteriors of alternating $3$-component pretzel links.

\begin{figure}
    \centering
     \includegraphics[width=0.7\linewidth]{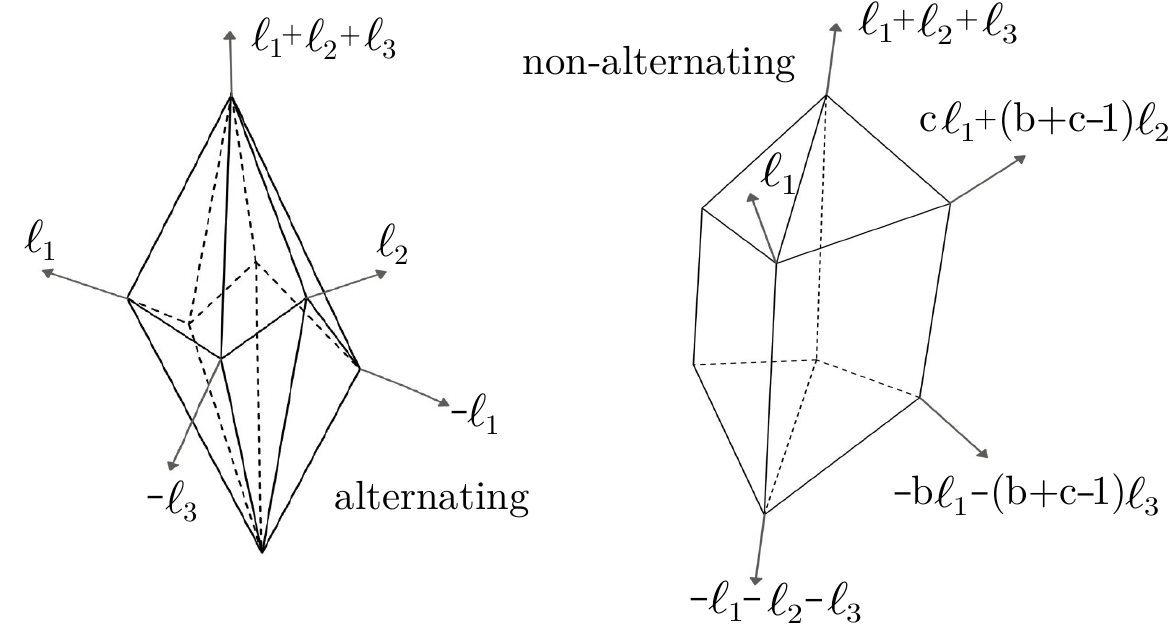}
     \caption{Left: the shape of the Thurston polyhedron for alternating pretzel links $P(2a,2b,2c)$. Right: the shape of the Thurston polyhedron for Pretzel links $P(2a,2b,2c)$ with $a<-1$, $b,c>0$ and $ab+bc+ac\le a$.}
     \label{fig: polyhedra}
\end{figure}

\bt\label{thm: alternating pretzels} Let $L\subset S^3$ be the pretzel link $P(2a,2b,2c)$, for some positive integers $a,b,c$. There is a choice of orientation and labelling for the components $\ell_1,\ell_2$ and $\ell_3$ of $L$ such that the Thurston ball of $H_2(M_L,\partial M_L;\R)\cong (\R\ell_1)\oplus(\R\ell_2)\oplus(\R\ell_3)$ is the polyhedron spanned by $$\pm \frac {\ell_1}{b+c-1},\pm \frac {\ell_2}{a+c-1},\pm \frac {\ell_3}{a+b-1}, \text{ and } \pm (\ell_1+\ell_2+\ell_3).$$
This polyhedron is depicted on the left of Figure \ref{fig: polyhedra}.
\et

By further analyses, we also determine the Thurston norm for some non-alternating pretzel links. 

\bigskip

\bt\label{thm: nonalternating pretzels} Let $L\subset S^3$ be the pretzel link $P(2a,2b,2c)$, for some integers $a<-1$ and $b,c>0$ satisfying $ab+bc+ac\le a$. There is a choice of orientation and labelling for the components $\ell_1,\ell_2$ and $\ell_3$ of $L$ such that the Thurston ball of $H_2(M_L,\partial M_L;\R)\cong (\R\ell_1)\oplus(\R\ell_2)\oplus(\R\ell_3)$ is the polyhedron spanned by $$\pm \frac {\ell_1}{b+c-1},\pm \frac {c\ell_1+(b+c-1)\ell_2}{(|a|-1)(b+c-1)},\pm \frac {b\ell_1+(b+c-1)\ell_3}{(|a|-1)(b+c-1)}, \text{ and } \pm (\ell_1+\ell_2+\ell_3).$$ 
This polyhedron is depicted on the right of Figure \ref{fig: polyhedra}.\et

Theorem \ref{thm: nonalternating pretzels} allows us to answer a question posed by Baker and Taylor in \cite{10.4310/jdg/1563242469}. Given a compact oriented $3$-manifold $M$ and a knot $K\subset M$, the \emph{wrapping number} with respect to $K$ is the map $$\wrap_K:H_2(M,\partial M)\to \Z_{\ge 0}$$ that takes a class $\alpha$ to the minimum number of intersections between $K$ and any norm-minimising incompressible surface representative of $\alpha$. In \cite[Question 2.1]{10.4310/jdg/1563242469}, Baker and Taylor asked whether the function $\wrap_K$ generally satisfies the triangle inequality. We can now answer this question negatively.

\vspace{8pt}

\noindent\textbf{Theorem \ref{thm: wrapping number}.} \textit{Let $L_0\subset S^3$ be a torus link $T(2,2c)$, for some $c>1$, and let $M_0:=M_{L_0}$ be its exterior. There exist infinitely many knots $K\subset M_0$ such that $\wrap_K$ does not satisfy the triangle inequality in $H_2(M_0,\partial M_0)$.}

\vspace{8pt}

Going back to the proposed strategy for computing the Thurston norm, one of its limitations is that a single hierarchy is often insufficient to fully recover the shape of the Thurston polyhedron. At the same time, finding hierarchies is a lengthy process and strongly depends on the capability of visualising $M$ or on the properties we know about it. However, given a groomed sutured manifold hierarchy $\mc H_+$ for $(M,\partial M)$, we can easily find a different groomed sutured manifold hierarchy by modifying $\mc H_+$. Indeed, suppose that the last decomposition of $\mc H_+$ is along a disk $D$. We can consider the hierarchy $\mc H_-$ that differs from $\mc H_+$ just in the chosen orientation of $D$. How much do the Euler classes of the two obtained foliations differ from each other? This question is answered by the following lemma.

\vspace{8pt}

\noindent\textbf{Lemma \ref{lemma: swapping orientation}.} \textit{Let $\mathcal{F_+}$ be a taut foliation arising from the groomed sutured manifold hierarchy $\mc H_+$ given by $$(M,\partial M)\rightsquigarrow ...\rightsquigarrow (M_r,\gamma_r)\overset{D}{\rightsquigarrow}  \text{product sutured balls}$$
where $D$ is a non-separating decomposing disk. Let $\mc H_-$ be the groomed sutured manifold hierarchy coinciding with $\mc H_+$ until $(M_r, \gamma_r)$, yet ending with $$(M_r,\gamma_r)\overset{\overline D}{\rightsquigarrow} \text{product sutured balls.}$$ Let $\mc F_-$ be a taut foliation carried by $\mc H_-$. The difference $e(\mathcal F_+)-e(\mathcal{F_-)}$ is the Poincar\'e dual of $(2-|D\cap \gamma_r|)i_*\delta$, where $i_*\delta\in H_1(M)$ is the image under the inclusion map $i_*:H_1(M_r)\to H_1(M)$ of the curve $\delta$ dual to $D$ in $M_r$.}  

\vspace{8pt}

As a consequence of the work of Gabai \cite{gabai_foliations_1987-1}, we deduce the following corollary (compare with \cite{miller, 10.4310/jdg/1563242469}).

\vspace{8pt}

\noindent\textbf{Corollary \ref{cor: linking number}.} \textit{Suppose that $M$ is the exterior of a nonsplit link $L\subset S^3$. Let $S$ be a taut surface in $M$ such that there is a boundary torus $T$ disjoint from $S$, and call $\ell$ the component of $L$ corresponding to $T$. If $[S]$ is not a corner of the Thurston norm on $H_2(M,\partial M;\R)$, then $\ell$ has linking number zero with every other component of $L$ and, for every non-longitudinal slope $\alpha$ on $T$, the surface $S$ remains norm-minimising in the Dehn-filled manifold $M_T(\alpha)$. Moreover, if $M$ is $S_T$-atoroidal, then $S$ is a leaf of a taut foliation of $M_T(\alpha)$ transverse to the core of the filling, for any non-longitudinal slope $\alpha\subset T$.}

\bigskip

\noindent \textbf{Structure of the paper.} Sections \ref{sec: thurston norm} and \ref{sec: hierarchies} are devoted to recalling the definitions and some of the main properties related to the Thurston norm, sutured manifold hierarchies, and branched surfaces. In Section \ref{sec: maw dual graph}, we explore the relations between Euler classes of taut foliations and the Thurston norm. Next, we recollect the features of the maw dual graph in Section \ref{subsec: maw dual graph} and prove Corollary \ref{cor: linking number}.

In Section \ref{sec: pretzels}, we propose an approach to compute the Thurston norm of any given manifold. Then, we prove Theorems \ref{thm: alternating pretzels} and \ref{thm: nonalternating pretzels}. Finally, in Section \ref{sec: wrapping number} we employ the result on pretzel links to prove Theorem \ref{thm: wrapping number}. 

\vspace{8pt}

\noindent \textbf{Notation.} If $P$ is a manifold, $|P|$ is the number of connected components of $P$. If $P$ is oriented and connected, $\overline P$ indicates $P$ with reversed orientation. 
The symbol $M$ generically indicates a smooth, compact and oriented $3$-manifold. The symbol $M^{\circ}$ refers to the interior of $M$, namely $M-\partial M$.
If $L$ is a link in $S^3$, the link exterior $S^3-N^{\circ}(L)$ is referred to as $M_L$.

If $S$ is a properly embedded surface in $M$, we use $N(S)$ to refer to a closed small tubular neighbourhood $S\times[-1,1]$ of $S$ in $M$. We will convey $\partial N(S)=S\times\{\pm 1\}$ and $N^{\circ}(S):=N(S)-\partial N(S)=S\times (-1,1)$. The manifold obtained by \emph{cutting} $M$ along $S$ is $M-N^{\circ}(S)$.

We will understand integral coefficients for the homology: for example, we will refer to $H_2(M,\partial M;\Z)$ just as $H_2(M,\partial M)$.
 
\vspace{8pt}

\noindent \textbf{Acknowledgements.} I am thankful to my advisor, Mehdi Yazdi, for his constant and insightful guidance. I would also like to thank Kenneth Baker and Diego Santoro for stimulating conversations.

\section{The Thurston norm}\label{sec: thurston norm}
We now recall some of the basic definitions and properties of the Thurston norm. The main results of this section, as for the use we will make of them during this article, are Lemmas \ref{lemma: additivity} and \ref{lemma: maximal cones}. 

\bigskip

\noindent Every properly embedded oriented surface $S\subset M$ represents a class $[S]\in H_2(M,\partial M)$. Conversely, every class in $H_2(M,\partial M)$ is represented by a properly embedded oriented surface. 

A natural question to ask is whether the sum in $H_2(M,\partial M)$ corresponds to some operation on surfaces in $M$. The answer is affirmative: if $[S_1]=\alpha_1$ and $[S_2]=\alpha_2$, then $\alpha_1+\alpha_2$ is represented by the \emph{oriented cut and paste} of $S_1$ and $S_2$. Here is the definition.

\begin{figure}
    \centering
    \includegraphics[width=0.7\textwidth]{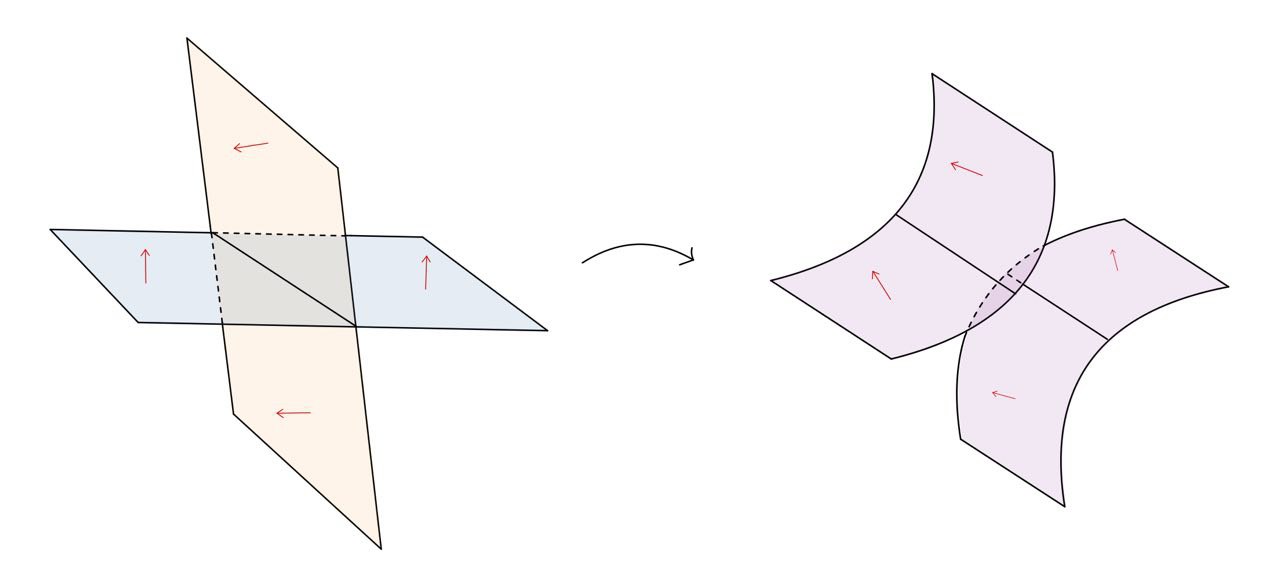}
    \caption{Local view of the oriented cut and paste operation. Red arrows indicate the normal orientation.}
    \label{fig: orientedcutandpaste}
\end{figure}

\bd(Oriented cut and paste) Let $S_1,S_2\subset M$ be properly embedded oriented surfaces that are transverse to each other. There is only one way to split $S_1$ and $S_2$ along $S_1\cap S_2$ and glue back to obtain an oriented surface $S$, which locally looks as in Figure \ref{fig: orientedcutandpaste}. The obtained surface $S$ is called the \emph{oriented cut and paste} of $S_1$ and $S_2$. Notice that $\chi(S)=\chi(S_1)+\chi(S_2)$.
\ed

\bd(Thurston norm) Let $S\subset M$ be a properly embedded orientable surface, and denote with $\chi(S)$ the Euler characteristic of $S$. Define $\chi_-(S):=|\chi(S')|$, where $S'$ is the surface obtained from $S$ after discarding all sphere or disk components.
For every class $\alpha\in H_2(M,\partial M)$, let $$x(\alpha)=\min_{[S]=\alpha}\chi_-(S),$$ where $S\subset M$ varies among the properly embedded oriented surfaces representing $\alpha$.
\ed

\bt[\cite{Thurston1986ANF}] The function $x:H_2(M,\partial M)\to \Z$ can be extended to a seminorm \\ $x:H_2(M,\partial M;\R)\to \R_{\ge 0}$, called the \emph{Thurston norm} associated to the manifold $M$.
\et 

\bd[Taut surfaces]\label{def: taut surfaces} A properly embedded oriented surface $S\subset M$ is \emph{norm-minimising} if $x([S])=\chi_-(S)$. The surface $S$ is taut if, furthermore, no subset of components of $S$ is homologically trivial, and any two parallel components of $\partial S$ are coherently oriented.
\ed

\br\label{rem: taut surfaces} Suppose $S$ is a norm-minimising surface in an irreducible and $\partial$-irreducible manifold $M$. Any compressible or $\partial$-compressible component of $S$ is null-homologous. Indeed, compressions and $\partial$-compressions increase the Euler characteristic of a surface without altering its homology class. Moreover, spheres, disks, compressible tori and ($\partial$-)compressible annuli are null-homologous because of the irreducibility assumptions. In particular, if $S$ is taut, then it is also incompressible and $\partial$-incompressible.
\er

\bexa Suppose $L\subset S^3$ is an oriented link with components $\ell_1,...,\ell_n$. By excision we get an identification $H_2(M_L,\partial M_L)\cong(\Z\ell_1)\oplus...\oplus(\Z\ell_n).$ Here, $[S]=a_1\ell_1+...+a_n\ell_n$ if and only if $a_i$ is the algebraic intersection number between $\partial S$ and the meridian of $\ell_i$, for each $i$.

    A surface representing $\ell_1+...+\ell_n$ is taut if and only if it is a minimal-genus Seifert surface for $L$.

    The class $\ell_1$ is represented by $S\cap M_L$, where $S\subset S^3$ is a Seifert surface for $\ell_1$. Therefore, computing $x(\ell_1)$ is analogous to computing the minimal genus of $\ell_1$ in $M_{\ell_2\cup...\cup \ell_n}$.
    \eexa

\bigskip

\subsection{Shape of the unit ball}

To compute a (semi)norm on $\R^n$ it is enough to know its unit ball. A useful feature of the Thurston norm is that its ball can be recovered by knowing a finite amount of data, thanks to the following result.

\bt[\cite{Thurston1986ANF}]\label{thm: ball is a polyhedron} The unit ball $B_x$ of the Thurston norm on $M$ is a finite (but possibly unbounded) convex polyhedron with rational vertices and symmetric with respect to the origin. 
\et

When $x$ is not a norm, i.e. when some connected surface with non-negative Euler characteristic is homologically non-trivial, the null-space of $x$ is a vector space contained in $B_x$. In particular, $B_x$ is bounded if and only if $x$ is a norm on $H_2(M,\partial M;\R)$. 

\bd[Boundary of the Thurston ball, cones and corners]
We refer to the set of points with Thurston norm one as the \emph{boundary} $\partial B_x$ of the Thurston ball. Since $B_x$ is a finite polyhedron, there are finitely many affine hyperplanes $H_1,...,H_k$ in $H_2(M,\partial M;\R)$ such that $$\partial B_x=\bigcup_{i=1}^k (\partial B_x\cap H_i).$$ Each nonempty intersection $\partial B_x\cap H_{i_1}\cap...\cap H_{i_j}$ is a \emph{closed face} of $\partial B_x$. An \emph{open face} is the interior of a closed face. If $F\subset \partial B_x$ is an open face, the \emph{open cone} over $F$ is the set $\{\lambda v\,|\; v\in F, \lambda\in \R_{ >0}\}$. If $F$ is instead a closed face, the \emph{closed cone} over $F$ is the closure of the cone over the interior $F$, defined as above. Following \cite{miller}, we say that a class $\alpha\in H_2(M,\partial M;\R)$ is a \emph{corner} of the Thurston norm if it does not lie in any top-dimensional open cone. Notice that each closed cone contains every class with vanishing Thurston norm. In particular, any such class is a corner of the Thurston norm. When the seminorm $x$ identically vanishes, we convey that there are exactly one closed cone and one open cone, both corresponding to $H_2(M,\partial M;\R)$.
\ed

A useful principle to investigate the shape of the Thurston ball of a manifold is the following.

\bl\label{lemma: additivity} Let $V$ be a finite-dimensional real vector space endowed with a seminorm $x$, whose unit ball $B_x$ is a finite (possibly unbounded) polyhedron. Two vectors $\alpha,\beta\in V$ lie in the closed cone over a common boundary face of $B_x$ if and only if $$x(\alpha+\beta)=x(\alpha)+x(\beta).$$
\el
\bp If $x(\alpha)=0$, both statements hold. So, we can suppose that neither $\alpha$ nor $\beta$ have zero seminorm. 

If $\alpha$ and $\beta$ lie in the closed cone over a common boundary face $F$ of $B_x$, then every point of the segment between $\alpha/x(\alpha)$ and $\beta/x(\beta)$ lies on $F$. In particular \begin{equation}\label{eq: additivity} x\left(\frac{x(\alpha)}{x(\alpha)+x(\beta)} \frac{\alpha}{x(\alpha)}+\frac{x(\beta)}{x(\alpha)+x(\beta)}\frac{\beta}{x(\beta)}\right)=1,\end{equation} and this is equivalent to $x(\alpha+\beta)=x(\alpha)+x(\beta)$.

Conversely, $x(\alpha+\beta)=x(\alpha)+x(\beta)$ implies the validity of equation (\ref{eq: additivity}). Thus, the three aligned points $$\frac{\alpha}{x(\alpha)},\frac{\beta}{x(\beta)} \text{ and } \frac{\alpha+\beta}{x(\alpha+\beta)}=\frac{x(\alpha)}{x(\alpha)+x(\beta)} \frac{\alpha}{x(\alpha)}+\frac{x(\beta)}{x(\alpha)+x(\beta)}\frac{\beta}{x(\beta)}$$ all lie on the boundary of the convex $B_x$, so the whole segment between $\alpha/x(\alpha)$ and $\beta/x(\beta)$ does.  
\ep

\bigskip

\subsection{Torus decompositions}

An incompressible torus $T$ in $M$ allows us to focus on the manifold obtained by cutting $M$ along $T$. It turns out that the Thurston norm behaves well with respect to this decomposition.

\bl\label{lemma: additive norm}\cite[Proposition 11.2]{lackenby_efficient_2021} Let $M$ be a compact oriented $3$-manifold with incompressible boundary. Let $T\subset M$ be a disjoint union of incompressible tori and let $M'$ be $M$ cut along $T$, seen as a submanifold of $M$. Consider a properly embedded oriented surface $S\subset M$ transverse to $T$. After indicating with $x_M$ and $x_{M'}$ the Thurston norms on $H_2(M,\partial M)$ and $H_2(M',\partial M')$ respectively, the following equality holds $$x_M([S])=x_{M'}([S\cap M']).$$
\el

When the torus $T$ is separating in $M$, the norm $x_{M'}$ is the sum of the norms of its connected components. In turn, the following lemma will apply.

\bl\label{lemma: maximal cones} Let $V, V_1$ and $V_2$ be finite-dimensional real vector spaces and $f_1:V\to V_1$ and $f_2:V\to V_2$ be linear maps. For $i=1,2$, let $x_i$ be a seminorm on $V_i$ whose unit ball $B_{x_i}$ is a finite (possibly unbounded) polyhedron. 

The function $x=x_1\circ f_1+x_2\circ f_2:V\to \R$ is a seminorm whose unit ball $B_x$ is a finite polyhedron. If $C_1$ and $C_2$ are closed cones over boundary faces of the unit balls of $x_1$ and $x_2$ respectively, then $f_1^{-1}(C_1)\cap f_2^{-1}(C_2)$ is a closed cone over a boundary face of the unit ball of $x$.
Vice versa, each closed cone $C$ over a boundary face of the unit ball of $x$ arises as an intersection $f_1^{-1}(C_1)\cap f_2^{-1}(C_2)$, for some $C_1$ and $C_2$ closed cones over boundary faces of the unit balls of $x_1$ and $x_2$ respectively.
%In particular, a class $\alpha\in V$ is in the interior of a top-dimensional cone of $x$ if and only if $f_i(\alpha)\in V_i$ is in the interior of a top-dimensional cone of $x_i$ for $i=1,2$.
\el
\bp The proof consists of a repeated use of Lemma \ref{lemma: additivity}. Given $\alpha,\beta\in V$: \begin{align*}    x(\alpha+\beta)&=x_1(f_1(\alpha+\beta))+x_2(f_2(\alpha+\beta))\\&=x_1(f_1(\alpha)+f_1(\beta))+x_2(f_2(\alpha)+f_2(\beta))\\ &\le x_1(f_1(\alpha))+x_1(f_1(\beta))+x_2(f_2(\alpha))+x_2(f_2(\beta))\\ &=x(\alpha)+x(\beta), 
\end{align*} and equality holds if and only if $f_1(\alpha),f_1(\beta)$ belong to a common closed cone of $x_1$ and $f_2(\alpha),f_2(\beta)$ belong to a common closed cone of $x_2$. Therefore, $\alpha,\beta$ belong to a common closed cone of $x$ if and only if $\alpha,\beta\in f_1^{-1}(C_1)\cap f_2^{-1}(C_2)$, for some closed cones $C_1,C_2$ of $x_1$ and $x_2$ respectively.
Since the closed cones over the faces of $x_i$ completely cover $V_i$ for $i=1,2$, the cones arising as intersections $f_1^{-1}(C_1)\cap f_2^{-1}(C_2)$ completely cover $V$, and $x$ is additive on each such intersection. This also shows that the unit ball of $x$ is a finite polyhedron. 

Suppose now that $C\subset V$ is a closed top-dimensional cone of $x$. For $i=1,2$, $x_i$ is additive on $f_i(C)$, hence there must be a closed cone $C_i$ of $x_i$ such that $f_i(C)\subset C_i$.
%This is because additivity on pairs (\alpha,\beta) guarantees additivity on tuples of any length within f_i(C). This in turn guarantees that any tuple belongs to a common closed cone by convexity. By the combinatorial structure of cones, one can find a finite generating set and use that as a tuple.
In particular, $C\subset f_1^{-1}(C_1)\cap f_2^{-1}(C_2)$. As $C$ is a maximal subset of $V$ on which $x$ restricts to a linear map, $C= f_1^{-1}(C_1)\cap f_2^{-1}(C_2)$. The lemma follows from the fact that every closed cone is an intersection of top-dimensional closed cones. 
\ep

We will use Lemma \ref{lemma: maximal cones} in the proofs of Proposition \ref{prop: linking number} and Corollary \ref{cor: linking number}.

\section{Taut sutured manifold hierarchies}\label{sec: hierarchies}
\subsection{Sutured manifolds}
Sutured manifolds were introduced by Gabai \cite{gabai_foliations_1984} to construct taut foliations on $3$-manifolds. In this section, we define the objects involved in the statements of Theorems \ref{thm: hierarchy} and \ref{thm: hierarchy to branched surface}. 

\begin{figure}
    \centering
    \includegraphics[width=0.8\linewidth]{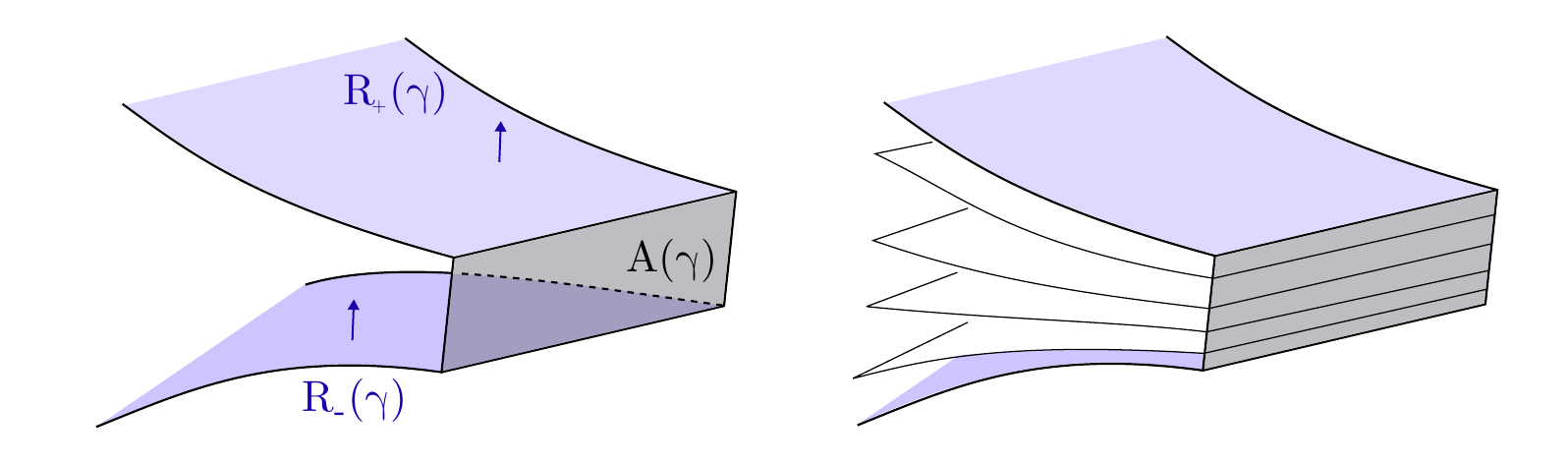}
    \caption{Left: local view of a sutured manifold close to an annular suture. Right: some leaves of a foliation, including the tangential boundary (in blue).}
    \label{fig: sutured manifold}
\end{figure}

Intuitively, a sutured manifold structure on a manifold $M$ is a decoration of $\partial M$ with convex corners, as in Figure \ref{fig: sutured manifold}. This corner structure is prescribing how a foliation is required to intersect $\partial M$, namely, transverse to $A(\gamma)$ and tangent elsewhere. Here is the rigorous definition of sutured manifolds.

\bd[Sutured manifold]\label{def: sutured manifold}  A \emph{sutured manifold} is a pair $(M,\gamma)$, where $M$ is a compact oriented $3$-manifold and $\gamma=A(\gamma)\sqcup T(\gamma)$ is a subsurface of $\partial M$ satisfying a number of conditions. The set $A(\gamma)$ is a union of pairwise disjoint annuli, whereas $T(\gamma)$ is a union of toroidal boundary components. A \emph{suture} is the core of an annulus of $A(\gamma)$, considered as an oriented simple closed curve. Every region of $R(\gamma):=\partial M-\gamma^{\circ}$ has a normal orientation. We call $R_+(\gamma)$ the union of the components of $R(\gamma)$ whose normal vector points out of $M$, and $R_-(\gamma)$ the union of the remaining components. The last requirement for $(M,\gamma)$ is that the normal orientations on $R(\gamma)$ are coherent with the orientation of the sutures: if $\delta\subset \partial R(\gamma)$ is a boundary component with the induced orientation, then $\delta$ is homologous to the respective suture in $H_1(\gamma)$.

We will say that $(M,\gamma)$ is a \emph{product sutured manifold}, if the sutured manifold $(M,\gamma)$ is homeomorphic to $(R\times[0,1],\partial R\times [0,1])$ for some compact oriented surface $R$ with non-empty boundary. 

If $\partial M$ is a union of tori and $S$ is a properly embedded oriented surface, call $(M(S),\gamma(S))$ the sutured manifold obtained by \emph{cutting} $(M,\partial M)$ along $S$; i.e., $M(S):=M-N^{\circ}(S)$ and $\gamma(S):=\partial M-N^{\circ}(S)$.
\ed

\bd[Taut foliations] A \emph{foliation} $\mathcal F$ of a sutured manifold $(M,\gamma)$ is a family of disjoint injectively immersed surfaces -- called \emph{leaves} of $\mathcal F$ -- that altogether cover $M$. We require leaves of $\mc F$ to be smoothly immersed in $M$ so that a small neighbourhood of each point is foliated as an open subset of $\R\times [0,+\infty)\times [-1,1],$ subdivided into horizontal half-planes $\R\times[0,+\infty)\times \{z\}$. Bearing this local model in mind, points in the interior of $\gamma$ correspond to points on $\R\times \{0\}\times(-1,1)$, points in the interior of $R(\gamma)$ correspond to points on $\R\times (0,+\infty)\times\{\pm 1\}$, and points on $\partial\gamma$ correspond to points on $\R\times\{0\}\times\{\pm 1\}$. 

The foliation $\mathcal F$ is \emph{cooriented} if $T\mc F$ is cooriented and the normal orientation points outward on $R_+(\gamma)$ and inward on $R_-(\gamma)$.

The foliation $\mc F$ is \emph{taut} if it is cooriented and every leaf intersects either a closed transversal or a properly embedded transverse arc with endpoints on $R(\gamma)$. %Gabai's theory produces taut foliations with an additional property, which we will hence include in the definition of taut foliations. We will require that the foliation $\mathcal F$ on $(M,\gamma)$ has no \emph{$2$-dimensional Reeb components} at the boundary, see Figure \ref{fig: 2-dim Reeb component}.
\ed

%\begin{figure}
%    \centering
%    \includegraphics[width=0.4\linewidth]{Figures/reeb_component.pdf}
%    \caption{A $2$-dimensional Reeb component is a subset of leaves of the boundary foliation that cover an annulus. The boundary of the annulus is composed of leaves with the same coorientation, i.e., they both point inside or both outside the annulus. The interior leaves are lines.}
%    \label{fig: 2-dim Reeb component}
%\end{figure}

The idea behind sutured manifolds is simple: take a properly embedded oriented surface $S$ in a manifold $M$ with toroidal boundary. We want to construct a foliation of $(M,\partial M)$ (i.e., a taut foliation transverse to the boundary) having $S$ as a union of leaves. Cut $M$ open along $S$ to obtain the sutured manifold $(M(S),\gamma(S))$. Here, both $R_+(\gamma(S))$ and $R_-(\gamma(S))$ are copies of $S$. Therefore, if we find a taut foliation of $(M(S),\gamma(S))$, we can then construct a taut foliation of $M$ with leaf $S$ by just glueing $R_+(\gamma(S))$ and $R_-(\gamma(S))$ back together.

Potentially, this technique can be iterated: say we find a new surface $S'$ in $M(S)$ and cut $M(S)$ open along $S'$. Then we define a new sutured manifold $(N',\gamma')$, where $N'=M(S)-N^{\circ}(S')$ and $\gamma'$ is suitably defined so that a foliation of $(N',\gamma')$ ``pulls back" to a foliation of $(M(S),\gamma(S))$. Next, we can decompose again $(N',\gamma')$ along a surface $S''$ and so on.

Finally, suppose that after some iterations we get a product sutured manifold $(R\times [0,1], \partial R\times [0,1])$. Then the process is done, because we can pull the product foliation back to a foliation of $(M,\partial M)$ having $S$ as a union of leaves.

The previous philosophy stands behind the notions of \emph{sutured manifold decompositions} and \emph{sutured manifold hierarchies}. See Section \ref{sec: pretzels} for concrete examples of this process.

The following definition is \cite[Definition 3.1]{gabai_foliations_1983} corrected as for \cite[Correction 0.3]{gabai_foliations_1987-1}.

\bd[Surface decompositions]\label{def: groomed} Let $(M,\gamma)$ be a sutured manifold. A \emph{decomposing surface} for $(M,\gamma)$ is a properly embedded oriented surface $S\subset M$ such that for every component $\lambda$ of $S\cap \gamma$ one of the following holds:
\begin{itemize}
    \item [$(i)$] $\lambda$ is a properly embedded nonseparating arc in $\gamma$;
    \item [$(ii)$] $\lambda$ is a simple closed curve in an annular component $A$ of $\gamma$ in the same homology class of the suture;
    \item [$(iii)$] $\lambda$ is a homotopically nontrivial curve in a toral component $T$ of $\gamma$, and if $\delta$ is another component of $T\cap S$, then $\lambda$ and $\delta$ represent the same homology class in $H_1(T)$.
\end{itemize}
We further require that:
\begin{itemize}
    \item [$(iv)$] No component of $\partial S$ bounds a disk in $R(\gamma)$ and no component of $S$ is a disk $D$ with $\partial D\subset R(\gamma)$.
\end{itemize}

\begin{figure}
    \centering
    \includegraphics[width=0.9\linewidth]{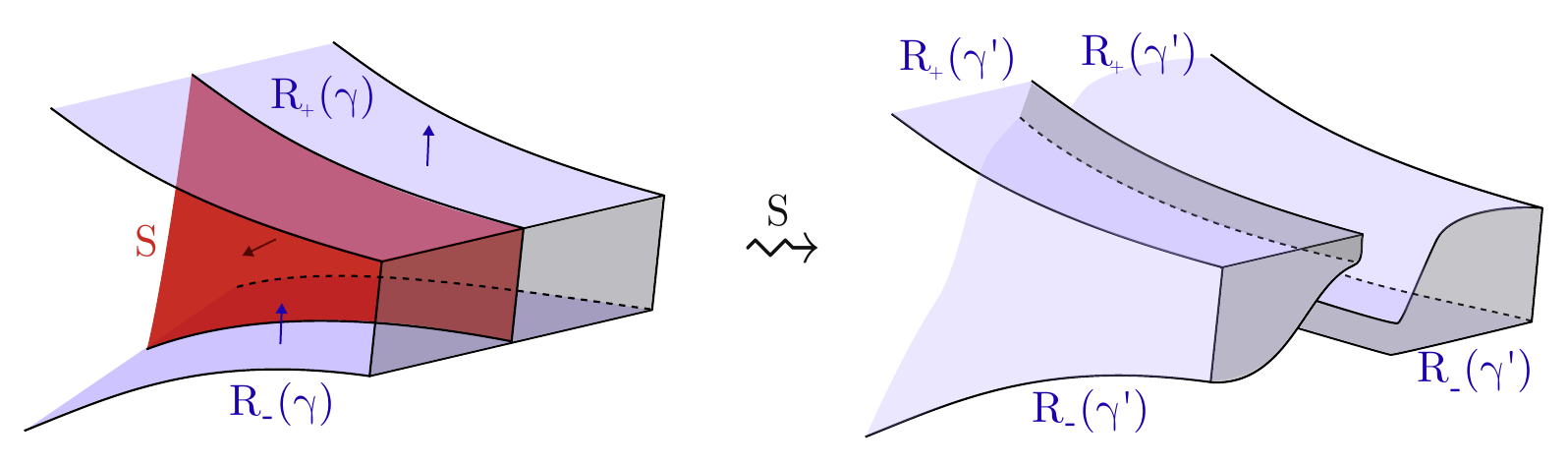}
    \caption{Local view of a surface decomposition close to a transverse intersection with a suture. On both left and right, the annular sutures $\gamma$ and $\gamma'$ are in grey.}
    \label{fig: surface decomposition}
\end{figure}

Any decomposing surface $S$ defines a \emph{sutured manifold decomposition} $(M,\gamma)\overset{S}{\rightsquigarrow}(M',\gamma')$, where $M'=M-N^\circ(S)$ and $R(\gamma')$ is obtained by smoothing together $\partial N(S)$ close to $R(\gamma)\cap M'$ according to the coorientation, as in Figure \ref{fig: surface decomposition}. More precisely: \begin{align*}
    \gamma'&=(\gamma\cap M')\cup N(S_+\cap R_-(\gamma))\cup N(S_-\cap R_+(\gamma)),\\
    R_+(\gamma')&= ((R_+(\gamma)\cap M')\cup S_+)-\overset{\circ}{N}(A(\gamma)),\\
    R_-(\gamma')&= ((R_-(\gamma)\cap M')\cup S_-)-\overset{\circ}{N}(A(\gamma)),
\end{align*}
where $S_{\pm}$ are the two components of $N(S)$, labelled according to coorientation. 

We say that the decomposing surface $S$ and the decomposition $(M,\gamma)\overset{S}{\rightsquigarrow}(M',\gamma')$ are \emph{groomed} if:
\begin{itemize}
    \item[$(v)$] No subset of toral components of $R(\gamma')$ is homologically trivial in $H_2(M)$, and
    \item[($vi$)] For each component $V$ of $R(\gamma)$, either:
    \begin{itemize}
        \item $S\cap V$ is a union of parallel, coherently oriented, nonseparating closed curves on $V$, or
        \item $S\cap V$ is a union of properly embedded arcs so that each component $\delta$ of $\partial V$ intersects $\partial S$ always with the same sign.
    \end{itemize}
\end{itemize} 

\ed

\bd[Hierarchies]\label{def: hierarchies} A \emph{taut sutured manifold hierarchy} for the sutured manifold $(M,\gamma)$ is a sequence of sutured manifold decompositions $$(M,\gamma)\overset {S_1}\rightsquigarrow...\overset{S_n}\rightsquigarrow (M_n,\gamma_n)$$ such that $(M_n,\gamma_n)$ is a union of product balls. The hierarchy is \emph{groomed} if each $S_i$ is groomed.  
\ed

We can now state a stronger version of Theorem \ref{thm: gabai intro}.

\bt\cite{gabai_foliations_1983}\label{thm: hierarchy} Let $M$ be a compact oriented and irreducible $3$-manifold $M$ with toroidal boundary and let $S\subset M$ be a properly embedded oriented surface. Suppose that parallel components of $\partial S$ are coherently oriented and no union of components of $S$ is zero in homology. The following are equivalent:
\begin{itemize}
    \item[(a)] $S$ is norm-minimising;
    \item [(b)] $(M,\partial M)$ admits a taut sutured manifold hierarchy starting with a decomposition along $S$ and so that no component of $R(\gamma_i)$ is a compressible torus;
    \item [(c)] $(M,\partial M)$ admits a groomed sutured manifold hierarchy starting with a decomposition along $S$ and so that no component of $R(\gamma_i)$ is a compressible torus;
    \item [(d)] $S$ is a union of leaves of a taut foliation of $(M,\partial M)$.
\end{itemize}
\et

Let us explain how Theorem \ref{thm: hierarchy} follows from \cite{gabai_foliations_1983}. Firstly, as already observed, $(M,\partial M)$ has a taut foliation with $S$ a union of leaves if and only if the cut-open manifold $(M(S),\gamma(S))$ admits a taut foliation. Consequently, the equivalence between conditions (a), (b) and (d) is established by \cite[Corollary 5.3]{gabai_foliations_1983} applied to the sutured manifold $(M(S),\gamma(S))$. Condition (c) clearly implies condition (b), whereas the converse follows from \cite[Lemma 5.4]{gabai_foliations_1983}.

%\bd A \emph{product disk} for $(M,\gamma)$ is a decomposing disk $D$ with $|\partial D\cap \gamma|=2$.\ed

%\br \begin{itemize}
 %   \item Every product sutured manifold can be decomposed into a disjoint union of product balls through a disjoint collection of product disks. In particular, every taut sutured manifold hierarchy can be completed, meaning that it can be extended to a complete hierarchy.
  %  \item The converse holds as well: if $(M,\gamma)$ is a sutured manifold with $R(\gamma)\ne \emptyset$ and it is completely decomposable through product disks, then $(M,\gamma)$ is a product sutured manifold \cite{gabai_foliations_1984}.
   % \item Decomposing along product disks is never a gamble: by \cite[Lemma 3.12]{gabai_foliations_1983}, in a product disk decomposition $(M,\gamma)\overset{D}{\rightsquigarrow}(M',\gamma')$, $(M,\gamma)$ admits a taut foliation if and only if $(M',\gamma')$ does.
%\end{itemize}
%\er

\subsection{Branched surfaces}

Branched surfaces are combinatorial tools to construct foliations in $3$-manifolds. We now quickly review their definition and how to reorganise the information of a sutured manifold hierarchy into a branched surface.

\begin{figure}
    \centering
    \includegraphics[width=\linewidth]{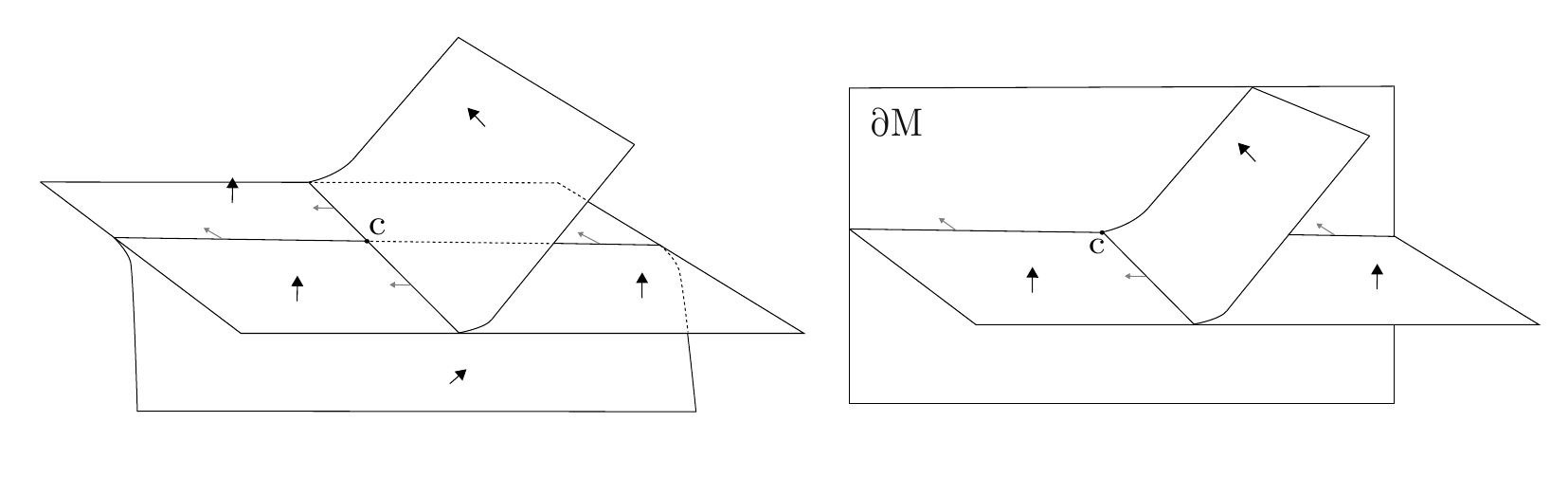}
    \caption{Local views of a branched surface close to a triple point $c$ (left) and at the intersection with $\partial M$ (right). Big arrows indicate local coherent choices of coorientations for the sectors, the small grey arrows indicate the maw vector field. In this case, the normal orientation of $\partial M$ is pointing out of $M$.}
    \label{fig: branching models}
\end{figure}

\bd[Branched surface] A \emph{branched surface} in a $3$-manifold $M$ is a $2$-complex $\mathcal B$ that locally looks like an open subset of the models in Figure \ref{fig: branching models}. Notice that a branched surface is naturally endowed with a tangent space at every point. Moreover, we always require branched surfaces to be transverse to $\partial M$.  

The \emph{branching locus} of $\mathcal{B}$ is the set $\brl(\mathcal B)$ of nonmanifold points of $\mathcal B$. It is a properly embedded graph on $\mathcal B$, whose vertices have degree either one or four, and the degree-four vertices are called \emph{triple points}. The branching locus $\brl(\mathcal B)$ is naturally endowed with the \emph{maw vector field}: this is any nonvanishing smooth vector field tangent to $\mathcal B$ and always pointing towards the branching direction. See again Figure \ref{fig: branching models}. Once a normal orientation on $\partial M$ is chosen, we extend the maw vector field on $\partial \mc B$ by requiring it to always point in the normal direction.

A \emph{sector} of $\mathcal B$ is the metric closure of any connected component of $\mathcal B -\brl(\mathcal B)$.

The branched surface $\mathcal B$ is \emph{cooriented} if there is a choice of normal direction on each sector of $\mathcal B$ that is coherent at the branching locus. Equivalently, $\mathcal B$ is cooriented if its tangent bundle $T\mathcal B$ is.
\ed

\begin{figure}
    \centering
    \includegraphics[width=0.7\linewidth]{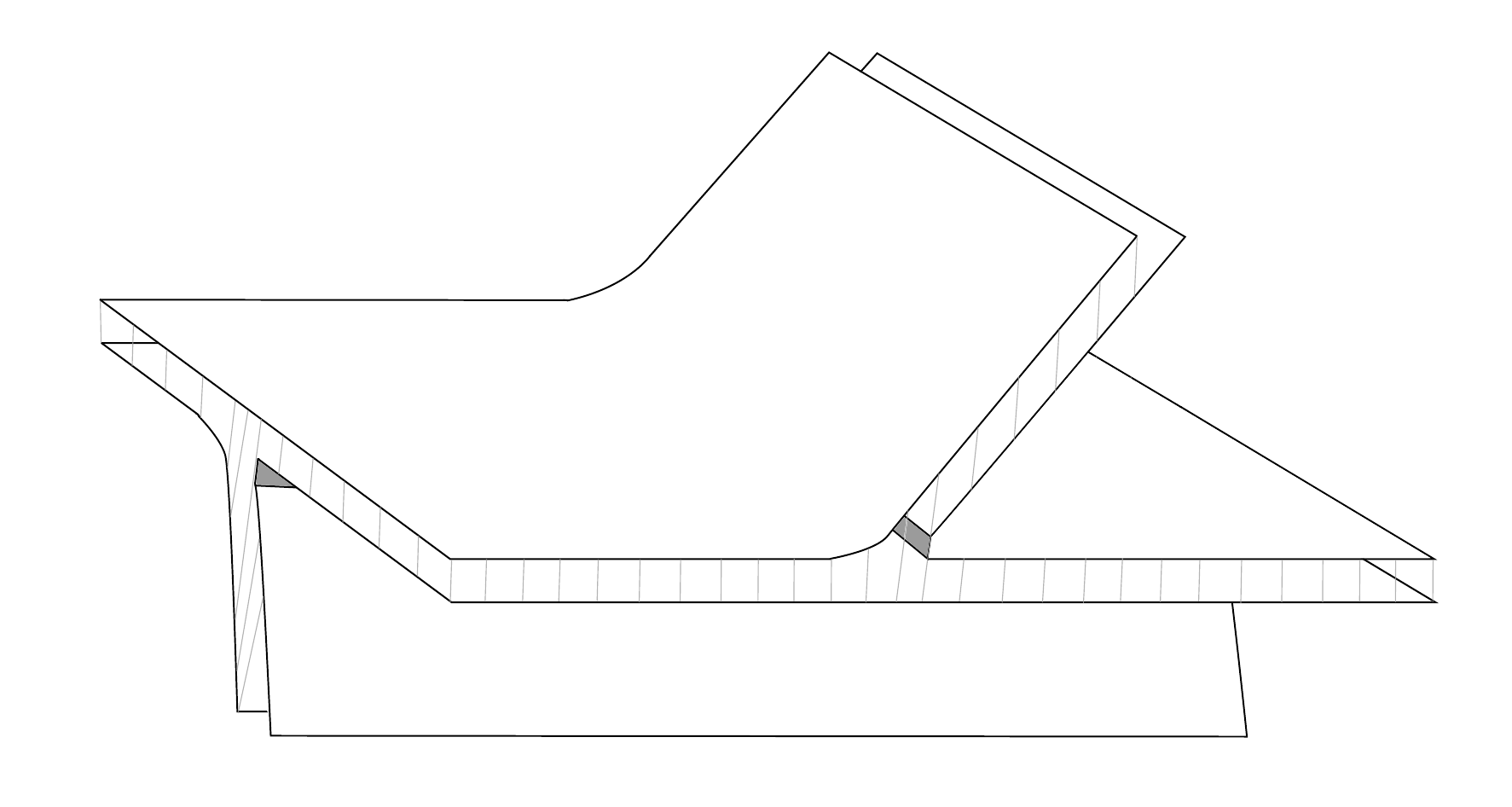}
    \caption{Local view of the fibered neighbourhood $N(\mc B)$ of a branched surface $\mc B$ close to a triple point of the branching locus. The grey intervals are some fibers. The small shaded areas are part of the vertical boundary $\partial_v N(\mc B)$.}
    \label{fig: fibered neighborhood}
\end{figure}

\bd[$I$-fibered neighbourhood and exterior] An \emph{$I$-fibered neighbourhood} of a branched surface $B\subset M$ is a regular neighborhood $N(B)$ foliated by interval fibers that intersect $B$ transversely, as in Figure \ref{fig: fibered neighborhood}. The boundary of $N(\mc B)$ has a \emph{horizontal} and a \emph{vertical} part. The \emph{horizontal boundary} $\partial_hN(\mc B)$ is given by the points on $\partial N(\mathcal B)$ that are transverse to the $I$-fibers. The \emph{vertical boundary} $\partial_vN(\mathcal B)$ consists of the portion of $\partial N(\mathcal B)\cap M^\circ$ that is tangent to the interval fibers.

The \emph{exterior} of $\mathcal B$ in $M$ is the sutured manifold $(M(\mathcal B),\gamma(\mathcal B))$ where:
\begin{itemize}
    \item $M(\mathcal{B})$ is the closure of $M-N(\mathcal B)$;
    \item $\gamma(\mathcal B)$ is the closure of $\partial M(\mc B)-\partial_hN(\mc B)$. In other words, $\gamma(\mc B)$ is the union of $\partial_vN(\mathcal B)$ and the closure of $\partial M- N(\mathcal B)$.
\end{itemize}
The connected components of $M(\mathcal B)$ are called \emph{complementary regions} of $\mathcal B$. 
\ed

\subsection{From sutured manifold hierarchies to branched surfaces}

Consider a sequence of sutured manifold decompositions $$(M,\partial M)\overset {S_1}\rightsquigarrow (M_1,\gamma_1)\overset {S_2}\rightsquigarrow...\overset{S_n}\rightsquigarrow (M_n,\gamma_n).$$ We can imagine each $M_i$ to be naturally embedded in $M$, and each surface $S_i\subset M$ with $\partial S_i=S_i\cap(\partial M\cup S_1\cup...\cup S_{i-1})$. For every $i=1,..,n$, we want to define a cooriented branched surface $\mc B_i\subset M$ so that $(M(\mc B_i), \gamma(\mc B_i))=(M_i,\gamma_i)$. This is the content of \cite[Construction 4.16]{gabai_foliations_1987}.

Set $\mc B_1=S_1$. Suppose then that $\mc B_i$ has been defined, and let us construct $\mc B_{i+1}$. After possibly a small isotopy, we can suppose $\partial S_{i+1}=S_{i+1}\cap (\partial M\cup \mc B_i)$ and that at most triple self-intersection points arise on $\mc B_i\cup S_{i+1}$. The coorientation on $\mc B_i$ and on $S_{i+1}$ gives instructions on how to smooth out $S_{i+1}$ close to $\mc B_i\cap S_{i+1}$ in order to get a cooriented branched surface $\mc B_{i+1}$. 

\begin{figure}
    \centering
    \includegraphics[width=0.7\linewidth]{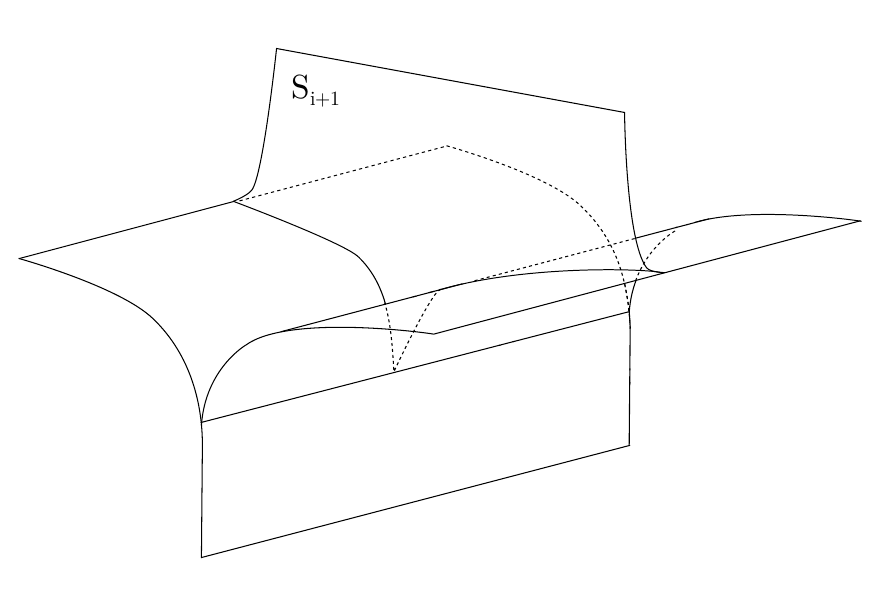}
    \caption{Local view of the branched surface $\mc B_{i+1}$ close to a transverse intersection between $S_{i+1}$ and $\partial _vN(\mc B_i)$.}
    \label{fig: third model}
\end{figure}

Notice that if $\partial S_{i+1}\cap \partial _vN(\mc B_i)\ne \emptyset$, the branched surface $\mc B_{i+1}$ will look like Figure \ref{fig: third model} near one such intersection. It requires a bit of thinking to see that this $2$-complex is indeed diffeomorphic to the left model in Figure \ref{fig: branching models} (which sectors are source, sink or none?).

%\newpage

\bt\cite[Construction 4.17]{gabai_foliations_1987}\label{thm: hierarchy to branched surface} If $$(M,\partial M)\overset {S_1}\rightsquigarrow (M_1,\gamma_1)\overset {S_2}\rightsquigarrow...\overset{S_n}\rightsquigarrow (M_n,\gamma_n)$$ is a groomed sutured manifold hierarchy, then the branched surface $\mc B_n\subset M$ fully carries a taut foliation. Namely, there is a taut foliation $\mc F$ of $M$ such that $\mc F \cap N(\mc B)$ is transverse to the interval fibers and $\mc F\cap M(\mc B)$ is the product foliation.
\et

\section{Euler classes of taut foliations and the Thurston norm}\label{sec: maw dual graph}

Let $M$ be a compact oriented manifold with (possibly empty) toroidal boundary. By \cite{novikov}, if $M$ has a taut foliation, then the components of $M$ which are not irreducible and $\partial$-irreducible are $S^2\times S^1$ or $D^2\times S^1$, foliated as products. We will henceforth assume that $M$ is irreducible and $\partial$-irreducible.

Suppose that $\mc F$ is a taut foliation of $(M,\partial M)$. Given a non-vanishing section $\sigma$ of the plane bundle $T\mathcal F|_{\partial M}\to \partial M$, the \emph{Euler class} of $\mathcal F$ relative to $\sigma$ is the obstruction $e(\mc F,\sigma)\in H^2(M,\partial M)$ to extending $\sigma$ to a non-vanishing section of the plane bundle $T\mc F\to M$. We refer the reader to \cite[Chapter 4]{candel_foliations_2003} for the rigorous definition of the Euler class, and to Corollary \ref{cor: maw dual graph} for an operative definition via an embedded graph in $M$. From now on, we will always understand $\sigma$ to be an outward-pointing section, and we will indicate $e(\mc F)=e(\mc F,\sigma)$.

In the present section, we establish Corollary \ref{cor: norm detection criterion}, which is one of the main ingredients in the proposed approach to computing the Thurston norm of a manifold. 

\subsection{Thurston's Inequality and its consequences}

As previously anticipated, the Euler class of a taut foliation and the Thurston norm on the ambient manifold are linked by Thurston's Inequality. 

\bt[Thurston's Inequality]\cite[Corollary 1]{Thurston1986ANF}\label{thm: thurston's inequality} Let $\mc F$ be a taut foliation of $(M,\partial M)$. For every $\alpha\in H_2(M,\partial M;\R)$, the following inequality holds
\begin{equation*}
    |\langle e(\mc F), \alpha\rangle| \le x(\alpha).
\end{equation*}
\et

\bd[Algebraically fully marked surfaces]\cite{gabai_fully_2020} Let $\mc F$ be a taut foliation of $(M,\partial M)$. Let $S\subset M$ be a properly embedded oriented surface with no sphere nor disk components. We say that $S$ is \emph{algebraically fully marked} with respect to $\mc F$ if the following equality holds \begin{equation*}
    \langle e(\mc F), [S]\rangle = \chi(S).
\end{equation*} When the foliation $\mc F$ is clear from the context, we just say that $S$ is algebraically fully marked.
\ed

\bc\cite{Thurston1986ANF}\label{cor: thurston} Let $\mc F$ be a taut foliation of $(M,\partial M)$. If $S\subset M$ is an algebraically fully marked surface, then $S$ is norm-minimising and $x([S])=-\langle e(\mc F),[S]\rangle$.
\ec
\bp Since $S$ has no sphere nor disk components, $\chi_-(S)=|\chi(S)|$. By Thurston's Inequality, $$\chi_-(S)=|\chi(S)|=|\langle e(\mc F), [S]\rangle|\le x([S]).$$ By the definition of the Thurston norm, we conclude that $\chi_-(S)=x([S])$.
\ep

\bc\label{cor: norm detection criterion} Let $\mathcal F$ be a taut foliation of $(M,\partial M)$. If $S$ and $S'\subset M$ are two algebraically fully marked surfaces, then their classes $[S]$ and $[S']$ lie on the same closed cone over a face of the Thurston ball.
\ec
\bp Suppose first that the oriented cut and paste $\hat S$ of $S$ and $S'$ has no sphere and no disk components. In this case, $\hat S$ is algebraically fully marked, because we have $$\chi(\hat S)=\chi(S)+\chi(S')=\langle e(\mc F), [S]\rangle+\langle e(\mc F), [S']\rangle=\langle e(\mc F), [\hat S]\rangle.$$ Since $$x(\hat S)=-\langle e(\mc F),[\hat S]\rangle=-\langle e(\mc F),[S]\rangle-\langle e(\mc F),[S']\rangle=x(S)+x(S'),$$ Lemma \ref{lemma: additivity} concludes the proof. We will now argue that we can always restrict to the case where $\hat S$ has neither sphere nor disk components. 

%Let $C$ be the annular and toroidal components of $S$. If $x([S])=0$ then the corollary automatically holds. Otherwise, notice that the classes $[S]$ and $[S-C]$ always belong to the same open cone over a face of the Thurston ball. We can therefore assume that all components of $S$ and $S'$ have negative Euler characteristic. Since compressions and $\partial$-compressions increase the Euler characteristic of a surface without altering its homology class, $S$ and $S'$ are incompressible and $\partial$-incompressible.
After possibly discarding null-homologous components, assume that $S$ and $S'$ are incompressible and $\partial$-incompressible (see Remark \ref{rem: taut surfaces}).
Suppose now that $\hat S$ has some sphere or disk component. An innermost argument shows that a component of $S\cap S'$ bounds a disk or is boundary-parallel in $S$ or $S'$. Since $M$ is irreducible and $\partial$-irreducible and $S,S'$ are incompressible and $\partial$-incompressible, we can remove such an inessential intersection via an isotopy of $S'$. After finitely many such isotopies, the oriented cut and paste of $S$ with the isotoped $S'$ has no sphere and no disk components.
\ep

\subsection{Dual Thurston norm and realisability of Euler classes}

After identifying $H^2(M,\partial M;\R)$ with $ \Hom (H_2(M,\partial M),\R)$, the Thurston norm $x$ on $H_2(M,\partial M;\R)$ induces a \emph{dual norm} $x^*$ defined by $$x^*(\beta):=\sup_{x(\alpha)\le 1} \langle \beta,\alpha\rangle $$ for every $\beta \in H^2(M,\partial M;\R)$, where $\alpha\in H_2(M,\partial M;\R)$ ranges among classes in the Thurston ball $B_x$. When $x$ is an honest norm, so is $x^*$. The dual unit ball $B_{x^*}$ is the polytope dual to the polyhedron $B_x$. Each open face $F$ of $B_x$ determines a dual face $F^*$ of $B_{x^*}$ with the property that equality $\langle \beta,\alpha\rangle=-x(\alpha)$ holds for every $\alpha$ in the closure of $F$ and for every $\beta$ in the closure of $F^*$.

Thurston's Inequality establishes that the Euler class of a taut foliation has dual Thurston norm at most one. Furthermore, it has norm one if the foliation has a compact leaf of negative Euler characteristic. Gabai showed that the converse is true for vertices of the dual Thurston ball.
The following theorem makes precise the statement that finitely many taut sutured manifold hierarchies are enough to recover the Thurston norm of a manifold.

\bt\cite{Gabai1997, gabai_fully_2020}\label{thm: gabai vertex} Let $M$ be a compact oriented irreducible $3$-manifold, possibly with toroidal boundary, and let $\beta\in H^2(M,\partial M;\R)$ be a vertex of the dual unit ball. Let $S$ be a taut surface in $M$ representing a class in the top-dimensional open cone dual to $\beta$. If $\mc F$ is a taut foliation such that $S$ is algebraically fully marked with respect to it, then $e(\mc F)=\beta$. In particular, $\beta$ is the relative Euler class of any taut foliation carried by a groomed sutured manifold hierarchy for $(M,\partial M)$ starting with a decomposition along $S$.
\et

\br\label{rem: vertex} The validity of Theorem \ref{thm: gabai vertex} can be extended to classes in the \emph{closed} top-dimensional cone dual to $\beta$. Namely, given a vertex $\beta\in H^2(M,\partial M;\R)$ of the dual unit ball and a taut surface $S$ in the closed top-dimensional cone dual to $\beta$, there is a groomed sutured manifold hierarchy for $(M,\partial M)$, starting with a decomposition along $S$, such that any carried foliation has relative Euler class $\beta$. 
This is because one can arrange that a groomed hierarchy starting with $S$ also carries a surface whose homology class lies in
the interior of the top-dimensional Thurston cone dual to $\beta$. Indeed, if $S'$ is any taut surface representing a class in the interior of the cone dual to $\beta$, then a sufficiently large oriented cut and paste $S'+nS$ yields a groomed decomposition for $(M(S),\gamma(S))$ by
\cite[Theorems 2.5 and 2.6]{Scharlemann1989}.
Consequently, the associated branched surface carries a surface in the interior
of the cone dual to $\beta$, and Gabai's uniqueness theorem implies that any
carried foliation has relative Euler class $\beta$.
\er
\section{The maw dual graph construction}\label{subsec: maw dual graph}
In the present section, we review the maw dual graph construction and prove Lemma \ref{lemma: swapping orientation}, as well as Corollary \ref{cor: linking number}.

Recall that $M$ is an irreducible and $\partial$-irreducible compact oriented $3$-manifold whose boundary is a (possibly empty) union of tori. Suppose $\mc F$ is a taut foliation carried by a cooriented branched surface $\mc B$ whose exterior is a union of balls (e.g., if $\mc B$ is constructed via a groomed sutured manifold hierarchy for $(M,\partial M)$). The \emph{maw dual graph} $\Gamma_m(\mc B)$ is an embedded simplicial $1$-cycle in $M$ representing the Poincar\'e dual of the relative Euler class $e(\mc F)\in H^2(M,\partial M)$. It can be constructed by suitably weighting the edges of the graph dual to the branched surface $\mc B$. We now recall this construction starting from the system of weights: the maw Euler characteristics of the sectors.

\begin{figure}
    \centering
    \includegraphics[width=0.8\linewidth]{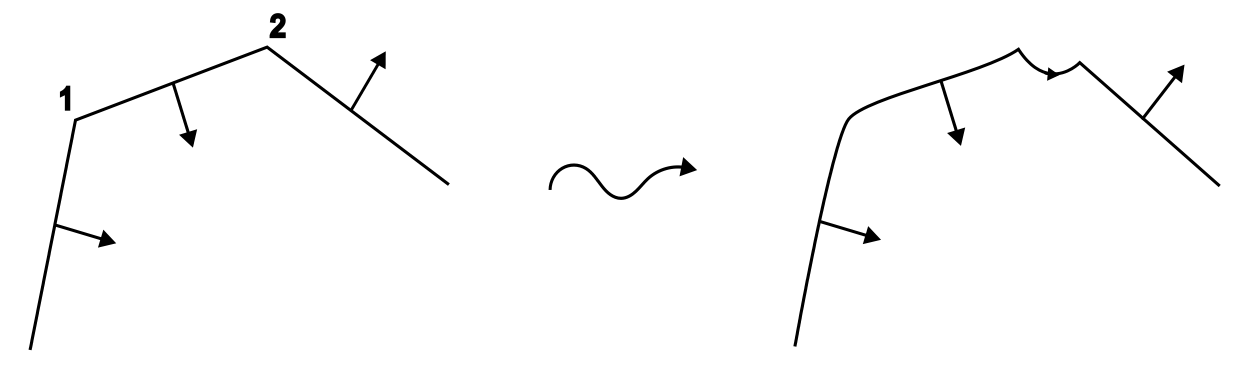}
    \caption{Left: a smoothable vertex, labelled (1), and a double corner, labelled (2). Right: the maw vector field induces a new corner structure on $\partial s$, where every smoothable vertex can indeed be smoothed (therefore inducing no corners), and each double vertex gives rise to a couple of right angles.}
    \label{fig: corners}
\end{figure}

\bd[Maw Euler characteristic] Let $s$ be a sector of a branched surface $\mathcal B$ and let $v$ be a vertex on $\partial s$, i.e., $v$ is either a triple point on $\partial s$ or an extreme of $s\cap \partial M$. We say that $v$ is a \emph{double corner} (or double vertex) for $s$ if the maw vector field points inward on one edge at $v$ and outward on the other. We say that $v$ is a \emph{smoothable corner} (or smoothable vertex) otherwise. See Figure~\ref{fig: corners}. Indicate with $\dc(s)$ the number of double corners on $\partial s$. The \emph{maw Euler characteristic} of $s$ is $$\chi_m(s):=\chi(s)-\frac 12 \dc(s).$$
\ed

\bd[Directed and maw dual graph]\label{def: maw dual graph} Let $\mathcal B$ be a cooriented branched surface in an ambient $3$-manifold $N$ with complement a union of balls. Fix a base point on each sector and inside each complementary region of $\mathcal B$. The \emph{directed dual graph} of $\mathcal B$ is the oriented graph $\Gamma(\mc B)\subset N$ constructed as follows:
\begin{itemize}
    \item Vertices: the vertices of $\Gamma(\mc B)$ are the base points in the complementary regions of $\mathcal B$;
    \item Edges: for every sector $s$ of $\mathcal{B}$, there is an arc $a(s)$ joining the vertices corresponding to the complementary regions of $\mathcal{B}$ adjacent at $s$. The arc $a(s)$ intersects $\mc B$ transversely and only at the base point of $s$. Moreover, $a(s)$ is oriented as the coorientation of $s$.
\end{itemize}

The \emph{maw dual graph} of $\mc B$ is the graph $\Gamma_m(\mc B)$ obtained by adding a system of weights to the edges of $\Gamma(\mc B)$. Namely, for each sector $s$ of $\mc B$, the dual arc $a(s)$ carries a weight equal to $\chi_m(s)$
\ed

\bt\cite[Theorem 1.2]{maw_dual_graph} Let $\mathcal{B}\subset N$ be a cooriented branched surface whose exterior is a union of product balls. Let $\mathcal{F}$ be a foliation fully carried by $\mathcal{B}$. The maw dual graph $\Gamma_m(\mc B)$ represents in $H_1(N)$ the Poincar\'e dual of the relative Euler class $e(\mathcal{F})\in H^2(N,\partial N)$.
\et

In particular:

\bc\label{cor: maw dual graph} Let $$(M,\partial M)\overset {S_1}\rightsquigarrow...\overset{S_n}\rightsquigarrow (M_n,\gamma_n)$$ be a groomed sutured manifold hierarchy, and let $\mc B\subset M$ be the induced branched surface.  If $\mc F$ is a taut foliation carried by the hierarchy, then the Poincar\'e dual of the relative Euler class $e(\mc F)\in H^2(M,\partial M)$ is represented by the maw dual graph $\Gamma_m(\mc B)$.
\ec

Another useful feature of the maw dual graph construction is the following property. 

\bl\cite[Lemma 4.7]{maw_dual_graph} Let $\mathcal{B'}\subset M$ be a cooriented branched surface whose exterior is a sutured manifold $(N,\gamma)$. Let $(D,\partial D)\subset (N,\partial N)$ be an oriented disk. Consider the cooriented branched surface $\mc B_+\subset M$ obtained by adding $D$ to $\mc B'$ (i.e., by suitably smoothing $D$ on $\mc B'$ about its boundary). Similarly, define the cooriented branched surface $\mc B_-$ by adding $\overline D$ to $\mc B'$. If the complement of $\mc B_
+$ in $M$ is a union of balls, then the symplicial $1$-chain $\Gamma_m(\mc B_+)-\Gamma_m(\mc B_-)$ is homologous to the chain $(2-|D\cap \gamma|)a_{\mc B_+}(D)$, where $a_{\mc B_+}(D)$ is the oriented arc of $\Gamma_m(\mc B_+)$ dual to the sector $D$.
\el

As a corollary, we obtain:

\bl\label{lemma: swapping orientation} Let $\mathcal{F_+}$ be a taut foliation arising from the groomed sutured manifold hierarchy $\mc H_+$ given by $$(M,\partial M)\rightsquigarrow ...\rightsquigarrow (M_r,\gamma_r)\overset{D}{\rightsquigarrow}  \text{product sutured balls}$$
where $D$ is a non-separating decomposing disk. Let $\mc H_-$ be the groomed sutured manifold hierarchy coinciding with $\mc H_+$ until $(M_r, \gamma_r)$, yet ending with $$(M_r,\gamma_r)\overset{\overline D}{\rightsquigarrow} \text{product sutured balls.}$$ Let $\mc F_-$ be a taut foliation carried by $\mc H_-$. The difference $e(\mathcal F_+)-e(\mathcal{F_-)}$ is the Poincar\'e dual of $(2-|D\cap \gamma_r|)i_*\delta$, where $i_*\delta\in H_1(M)$ is the image under the inclusion map $i_*:H_1(M_r)\to H_1(M)$ of the curve $\delta$ dual to $D$ in $M_r$.
\el

In the previous lemma, the fact that the hierarchy $\mc H_-$ is taut comes from the fact that the component of $M_r$ containing $D$ must be a sutured solid torus. Moreover, $D$ must be a meridional disk and $\delta$ the core of the solid torus. The details are left to the reader as an exercise on sutured manifolds.

\subsection{Exceptional slopes}

Suppose now that there is a taut surface $S\subset M$ disjoint from a boundary torus $T$. Given a slope $\alpha$ on $T$, we might wonder whether $S$ remains norm-minimising in the manifold $M_T(\alpha)$ obtained by Dehn-filling $T$ along the slope $\alpha$. In \cite{gabai_foliations_1987-1}, Gabai addressed this problem. Namely, he showed that there is at most one slope $d$ such that $S$ is not norm-minimising in $M_T(d)$. If such a slope exists, we call it the \emph{exceptional slope} of $T$ with respect to $S$. Under the stronger assumption that $M$ is \emph{$S_T$-atoroidal} -- i.e., if every torus $I$-cobordant with $T$ in $M-S$ is in fact parallel to $T$ -- Gabai showed that $S$ is a leaf of a taut foliation of $M_T(\alpha)$ transverse to the core of the filling, for every non-exceptional $\alpha$. A closer look at his construction enables us to deduce the following result.

\bprop\label{prop: linking number} Let $S$ be a taut surface in $M$ disjoint from a boundary torus $T$. If $[S]$ is not a corner of the Thurston norm on $H_2(M,\partial M;\R)$, then the inclusion map $H_1(T;\R)\to H_1(M;\R)$ is not injective. Furthermore, if $d$ is exceptional for $S$, then it generates the kernel of $H_1(T;\R)\to H_1(M;\R)$.
\eprop

\bc\label{cor: linking number} Suppose that $M$ is the exterior of a nonsplit link $L\subset S^3$. Let $S$ be a taut surface in $M$ such that there is a boundary torus $T$ disjoint from $S$, and call $\ell$ the component of $L$ corresponding to $T$. If $[S]$ is not a corner of the Thurston norm on $H_2(M,\partial M;\R)$, then $\ell$ has linking number zero with every other component of $L$ and, for every non-longitudinal slope $\alpha$ on $T$, the surface $S$ remains norm-minimising in $M_T(\alpha)$. Moreover, if $M$ is $S_T$-atoroidal, then $S$ is a leaf of a taut foliation of $M_T(\alpha)$ transverse to the core of the filling, for every non-longitudinal slope $\alpha\subset T$.
\ec

 Notice that a surface $S$ as above arises when $M$ is the exterior of a nonsplit link in $S^3$, and two components of the link have zero linking number. Indeed, call $\ell_1$ and $\ell_2$ such components. Let $S$ be the intersection of a Seifert surface for $\ell_1$ with the link exterior $M_L$. Since $\ell_1$ and $\ell_2$ have zero linking number, $\partial S$ intersects the boundary component $T$ of $M_L$ relative to $\ell_2$ in a family of $k\ge 0$ positively oriented meridians and $k$ negatively oriented meridians. After glueing $k$ tubes to $S$ along these meridians, we get a surface $S'$ in $M_L$ homologous to $S$ and disjoint from the boundary torus $T$.

\bp[Proof of Proposition \ref{prop: linking number}] \underline{Case $1$: $M$ is $S_T$-atoroidal.} We construct two taut foliations $\mathcal F_+$ and $\mathcal F_-$ on $M$, having $S$ as a leaf and whose relative Euler classes differ by a nonzero multiple of the Poincar\'e dual of the candidate exceptional slope $d\subset T$. Since $[S]$ is in the interior of a top-dimensional cone, the classes $e(\mc F_\pm)$ must coincide with the dual vertex on $B_{x^*}$, hence $[d]$ vanishes in $H_1(M;\R)$. The fact that $[d]$ generates the kernel of the inclusion map comes from the general fact that this map cannot vanish identically. %this is because if d vanishes, then c (defined below) intersects a cycle with boundary d exactly once, hence c is nonzero in M

\begin{figure}
    \centering
    \includegraphics[width=0.4\linewidth]{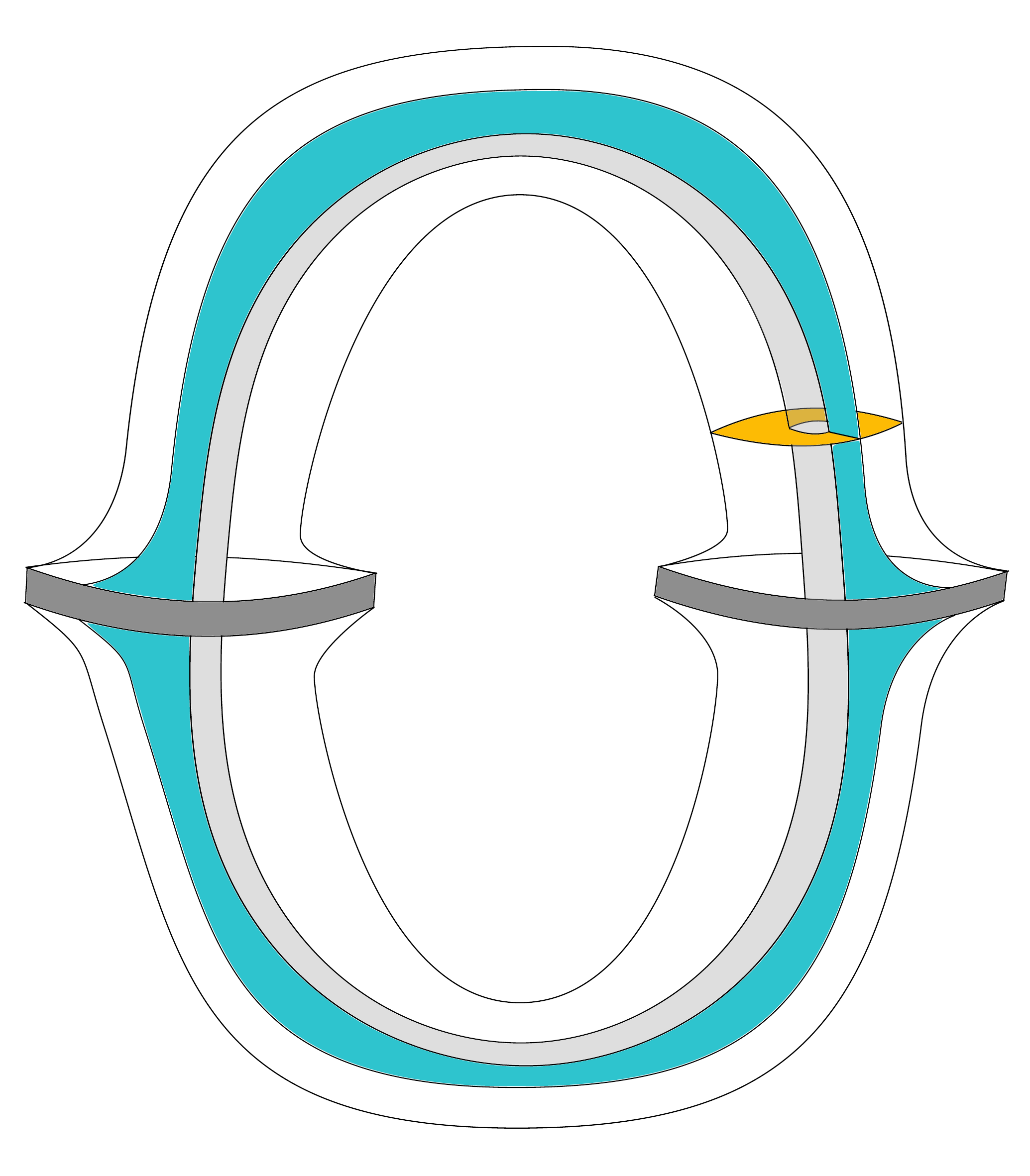}
    \caption{The annuli $A_d$ (orange) and $A_c$ (light blue), inside the sutured manifold $(H,\delta)$. To induce a groomed sutured decomposition, the outermost boundary of the orange annulus $A_d$ should be moved to one of the dark grey sutures. This is here avoided for a clearer picture.}
    \label{fig: linking number}
\end{figure}

As in \cite[Proof of Theorem 1.7]{gabai_foliations_1987-1}, there is a sequence of sutured manifold decompositions \begin{equation}\label{gabai's hierarchy}
    (M,\partial M)\overset{S}{\rightsquigarrow}(M(S),\gamma(S))\rightsquigarrow...\rightsquigarrow(M_k,\gamma_k)
\end{equation} so that $(M_k,\gamma_k)$ is a disjoint union of product balls and a sutured manifold $(H,\delta)$. The manifold $H$ is a collar $T\times [0,1]$ of $T$ in $M$, with $T=T\times 0$ and $\emptyset\ne\delta\subset T\times 1$. The sutures in $\delta$ are all isotopic to a simple closed curve $d\subset T\times 1$ and any exceptional slope for $S$ must coincide with $d$. By \cite[Lemma 5.4]{gabai_foliations_1983}, it is not restrictive to assume that the sequence (\ref{gabai's hierarchy}) is groomed. Let $c\subset T\times 1$ be a simple closed curve intersecting $d$ once. Consider the two annuli $A_d=d\times[0,1]\subset H$ and $A_c=c\times[0,1]\subset H$. See Figure \ref{fig: linking number}.

Choose an orientation on $A_d$ and decompose $(H,\delta)$ along it to obtain a sutured solid torus $T_d$. The annulus $A_c$ restricts to a meridional disk for $T_d$ intersecting the sutures $|\delta|+2$ times. Both choices of orientation on $A_c$ carry taut foliations $\mc F_{\pm}$ on $M$, and the surface $S$ is a union of leaves of these foliations. By Lemma \ref{lemma: swapping orientation}, the Euler classes of $\mc F_\pm$ differ by the Poincar\'e dual of $|\delta|[d]\in H_1(M)$. Since $[S]$ is not a corner of the Thurston norm, $[d]$ vanishes in $H_1(M;\R)$.

\bigskip

\underline{Case $2$: general case.} There is a torus $T'\subset M-S$ that separates $M$ into two pieces $M'$ and $V$, such that $S\subset M', T\subset V$, $M'$ is $S_{T'}$-atoroidal and $V$ is an $I$-cobordism between $T$ and $T'$. Since $M$ is irreducible, the torus $T'$ is incompressible: if it was compressible, it would either bound a solid torus or be contained in a $3$-ball. The existence of $T$ and $S$ rules out both possibilities. Lemma \ref{lemma: additive norm} now guarantees that for every $\alpha\in H^1(M;\R)$ the following formula holds $$x_M(\alpha)=x_{M'}(i_{M'}^*\alpha)+x_V(i_V^*\alpha),$$ where $x_P$ indicates the Thurston norm on $H^1(P;\R)\cong H_2(P,\partial P;\R)$, and $i_P^*:H^1(M;\R)\to H^1(P;\R)$ is the map induced by the inclusion $P\subset M$, for $P=M',V$.

By Lemma \ref{lemma: maximal cones}, each closed cone of $B_{x_M}$ arises as an intersection $(i_{M'}^*)^{-1}(C_{M'})\cap (i_{V}^*)^{-1}(C_{V})$, with $C_{M'}$ and $C_V$ closed cones of $x_{M'}$ and $x_V$ respectively. Moreover, the map $i_{M'}^*$ is open and surjective because $V$ is an $I$-cobordism. As the class $[S]$ is not a corner of $x_{M}$, then $[S]$ is not a corner of $x_{M'}$ either. By \cite[Proof of Corollary 2.4]{gabai_foliations_1987-1}, if $d$ is the exceptional slope of $T$ with respect to $S$ in $M$, then $d$ is homologous in $H_2(V,\partial V;\R)$ to the exceptional slope of $T'$ with respect to $S$ in $M'$. The previous case concludes the proof. \ep

\bexa Let $S$ be a minimal-genus Seifert surface for a knot $K_0\subset S^3$, and consider a taut foliation $\mc F_0$ of the exterior $M_0$, having $S$ as a leaf. Let $K\subset S$ be a separating essential simple closed curve, let $S_1,S_2$ be the components of $S-K$ in the drilled $3$-manifold $M:=M_0-N^\circ(K)$. In other words, $M$ is the exterior in $S^3$ of the link $K_0\cup K$. The surface $S$ can be pushed out of $N(K)$, hence embedded in $M$. The resulting taut representative of $[S]\in H_2(M,\partial M)$ is homologous to $S_1\cup S_2$, hence not a corner of the Thurston norm. In this case, one can readily construct taut foliations of the Dehn-filled $M_{\partial N(K)}(\alpha)$ having $S$ as a leaf, as long as $\alpha$ is not the longitude of $K$. Namely, first modify $\mc F_0$ by substituting the leaf $S$ with a product-foliated neighbourhood $S\times [-1,1]$. Then, excise a small annular neighbourhood $K\times (-1,1)\subset S\times 0$ from $S\times 0$ and blow air in it to obtain a foliation $\mc F$ of $M$. By means of classical operations on foliations of $3$-manifolds, a taut foliation of $M_{\partial N(K)}(\alpha)$ can be obtained by stacking saddles in the filling solid torus and performing $I$-bundle replacement along curves in $S_1\cup S_2$.
\eexa

\section{An approach for computing the Thurston norm}\label{sec: pretzels}
Thurston's Inequality makes Corollary \ref{cor: maw dual graph} relevant for getting information on the Thurston unit ball of a manifold. In what follows, consider the case of a compact oriented irreducible $3$-manifold $M$ with incompressible toroidal boundary (e.g. the exterior of a nonsplit link in the $3$-sphere).

\bigskip

\noindent Here is an approach for computing the Thurston norm of $H_2(M,\partial M;\R)$ by means of sutured manifold hierarchies:

\begin{enumerate}
    \item Start by a surface $S$ that we believe should be taut. 
    \item If $S$ is indeed taut, we can find a groomed sutured manifold hierarchy $\mc H$ for $(M,\partial M)$ starting with a decomposition along $S$, thanks to Theorem \ref{thm: hierarchy}.
    \item Reorganise the data of $\mc H$ into a branched surface $\mc B$ and compute the homology class of the maw dual graph $[\Gamma_m(\mc B)]\in H_1(M)$. The Thurston unit ball of $H_2(M,\partial M;\R)$ is contained in the strip $\{\alpha\in H_2(M,\partial M;\R)\,|\; |\langle \alpha,[\Gamma_m(\mc B)]\rangle|\le 1\}$. This is because of Thurston's Inequality (Theorem \ref{thm: thurston's inequality}) and Corollary \ref{cor: maw dual graph}.
    \item Focus on other classes $\alpha_1,...,\alpha_k\in H_2(M,\partial M)$ for which we know surface representatives. Thurston's Inequality gives a lower bound  $x(\alpha_i)\ge |\langle \alpha_i,[\Gamma_m(\mc B)]\rangle|$ for each $i=1,..,k$. If we manage to find surfaces $S_1,...,S_k$ representing $\alpha_1,...,\alpha_k$ respectively and each satisfying $\chi([S_i])= \langle S_i,\Gamma_m(\mc B)\rangle$, then they are norm-minimising. Moreover, the positive cone spanned by the classes $[S],\alpha_1,...,\alpha_k$ in $H_2(M,\partial M;\R)$ is contained in a closed cone of the Thurston norm by Corollaries \ref{cor: norm detection criterion} and \ref{cor: maw dual graph}.  
    \item Iterate the process above by constructing new hierarchies starting with a decomposition along $S$ or $S_i$, for some $i$. For instance, a new sutured manifold hierarchy can be constructed by reversing the orientation of a decomposing disk in the last decomposition of $\mc H$. The Euler class of this modified hierarchy can then be easily computed by means of Lemma \ref{lemma: swapping orientation}. 
\end{enumerate}

\emph{Provided we know ``enough surfaces" in $M$}, thanks to Theorem \ref{thm: gabai vertex} and Remark \ref{rem: vertex}, the above strategy can be repeated until one finds enough cones to cover all $H_2(M,\partial M;\R)$, hence fully recovering the Thurston norm. 
In the case of link exteriors in $S^3$, we can often find many surfaces and try to run the method above. In the case of pretzel links with three components, we obtain:

\bt\label{thm: 3-pretzel norm} Let $L\subset S^3$ be the pretzel link $P(2a,2b,2c)$, for some nonzero integers $a,b,c$. Orient and label $\ell_1,\ell_2$ and $\ell_3$ the components of $L$ as in the left of Figure \ref{fig: pretzels}. After possibly mirroring and relabelling, assume $b,c> 0$. Let $B_x$ be the Thurston ball of $H_2(M_L,\partial M_L;\R)\cong (\R\ell_1)\oplus(\R\ell_2)\oplus(\R\ell_3)$. 
\begin{itemize}
    \item If $a>0$ (i.e., when the standard diagram of $L$ is alternating), then $B_x$ is the polyhedron spanned by $$\pm \frac {\ell_1}{b+c-1},\pm \frac {\ell_2}{a+c-1},\pm \frac {\ell_3}{a+b-1}, \text{ and } \pm (\ell_1+\ell_2+\ell_3).$$
    \item If $ab+bc+ac\le a\le -2$, then $B_x$ is the polyhedron spanned by $$\pm \frac {\ell_1}{b+c-1},\pm \frac {c\ell_1+(b+c-1)\ell_2}{(|a|-1)(b+c-1)},\pm \frac {b\ell_1+(b+c-1)\ell_3}{(|a|-1)(b+c-1)}, \text{ and } \pm (\ell_1+\ell_2+\ell_3).$$ 
\end{itemize}
\et

As observed by Mehdi Yazdi, the condition $ab+bc+ac\le a$ automatically holds if $-a\ge \min (b,c)$. This is because one can rewrite $ab+bc+ac\le a$ as $-a\ge \frac{bc}{b+c-1}$, and the RHS of this expression is bounded above by $\min(b,c)$.

The proof below illustrates a general method that reduces the argument to a long but routine computational exercise. The reader may wish to concentrate on the structure of the argument rather than on the individual computations.

\begin{figure}
    \centering
    \includegraphics[width=\linewidth]{Figures/pretzels.pdf}
    \caption{Left: the diagram $D$ of the pretzel link $L=P(2a,2b,2c)$. Right: the Seifert surface $S$ of the link $L$. Topologically, $S$ is a pair of pants.}
    \label{fig: pretzels}
\end{figure}

\bp Consider the diagram $D$ of $L$ in Figure \ref{fig: pretzels} left. The application of Seifert's algorithm to the diagram $D$ yields the Seifert surface $S$ in Figure \ref{fig: pretzels} right. A handle decomposition of $S$ consists of two $0$-handles and three $1$-handles corresponding to the twisting regions. Therefore $S$ has Euler characteristic $-1$ and should reasonably be taut as $-1$ is the maximal Euler characteristic for a compact oriented and connected surface with three boundary components.

\bigskip

\begin{figure}
    \centering
    \includegraphics[width=0.8\linewidth]{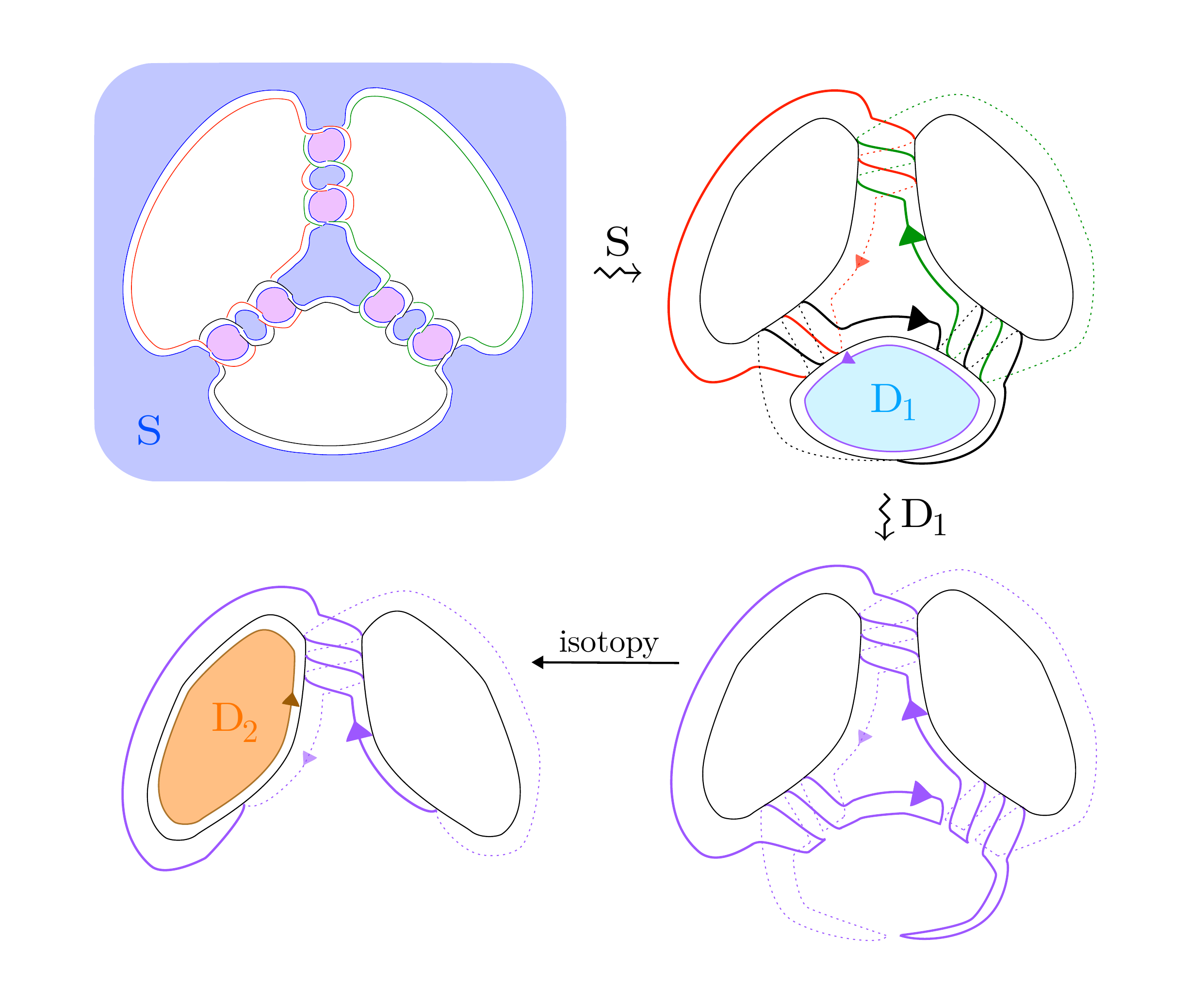}
    \caption{Top left: the Seifert surface $S$ embedded in the exterior of the link $L=P(-4,4,4)$. Top right: the sutured manifold $(M(S),\gamma(S))$ obtained by decomposing $(M_L,\partial M_L)$ along $S$. The decomposing disk $D_1$ is coloured light blue. The black, red and green curves represent the sutures of $\gamma(S)$. Bottom right: the sutured manifold $(M_1,\gamma_1)$ resulting by decomposing $(M(S),\gamma(S))$ along $D_1$. Bottom left: the manifold $(M_1,\gamma_1)$ after an isotopy of $\gamma_1$. The decomposing disk $D_2$ is coloured orange. The reader can verify that, after decomposing along $D_2$, the sutured manifold $(M_1,\gamma_1)$ yields a sutured manifold homeomorphic to a ball with exactly one suture, hence a product ball.}
    \label{fig: pretzel hierarchy}
\end{figure}

\noindent\underline{The first hierarchy.} Let us find a groomed sutured manifold hierarchy $\mc H_{--}$ for the link exterior $(M_L,\partial M_L)$ starting with the decomposition along $S$. The manifold $M(S)=M_L-N^\circ(S)$ is the complement in $S^3$ of a genus-$2$ handlebody, as in Figure \ref{fig: pretzel hierarchy} top right. Moreover, $M(S)$ is a genus-$2$ handle body as well, and two disk decompositions yield a manifold homeomorphic to a ball. This is illustrated in Figure \ref{fig: pretzel hierarchy}. After decomposing $(M(S),\gamma(S))$ along the disk $D_1$, we get a sutured solid torus $(M_1,\gamma_1)$. This can be further decomposed into a $3$-ball along a meridional disk $D_2$. As the resulting $3$-ball has only one suture, it is a product ball.

\begin{figure}
    \centering
    \includegraphics[width=0.8\linewidth]{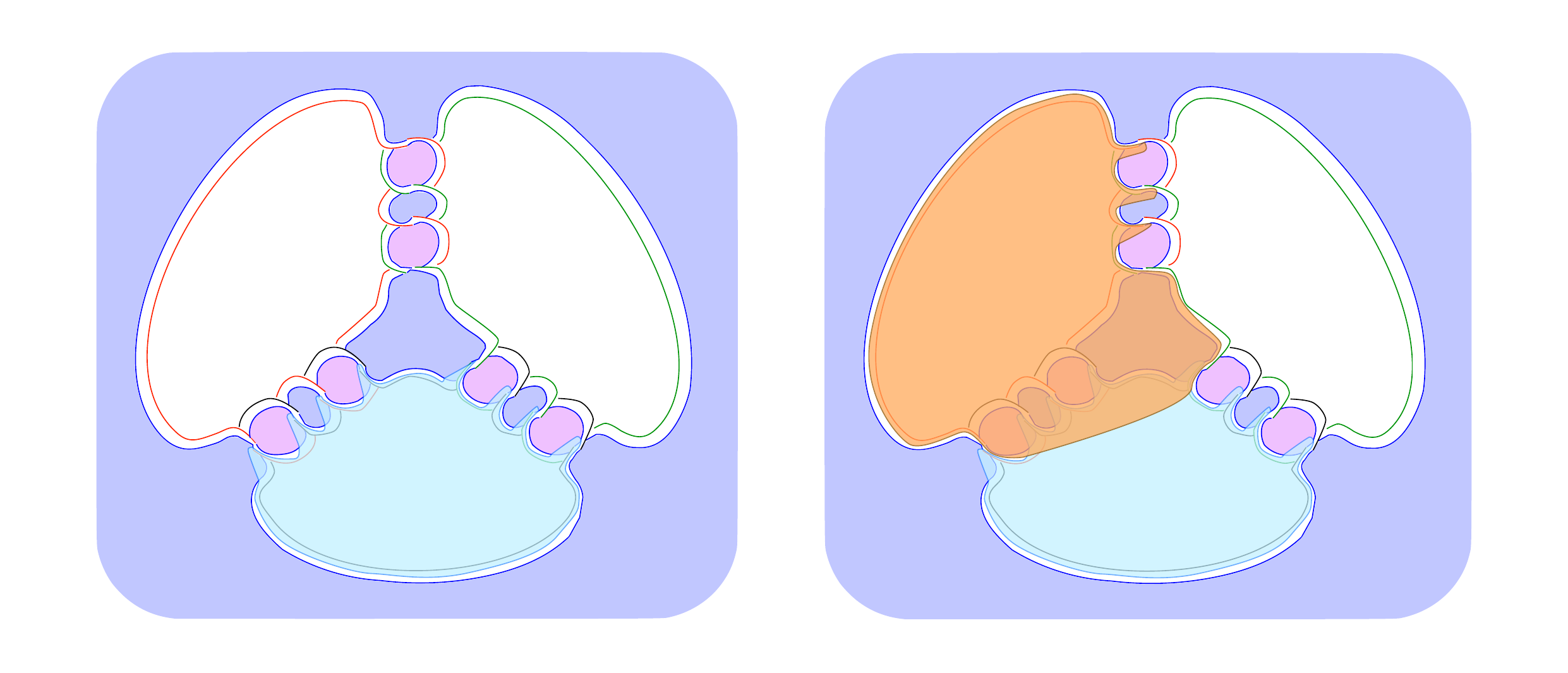}
    \caption{Left: the branched surface $\mc B'_{--}$. Right: the branched surface $\mc B_{--}$. The coorientations of $D_1, D_2$ and of the blue side of $S$ are pointing out of the page.}
    \label{fig: pretzel branched surface}
\end{figure}

Let $\mc B_{--}$ be the branched surface associated with the groomed sutured manifold hierarchy $\mc H_{--}$: $$(M_L,\partial M_L)\overset{S}{\rightsquigarrow} (M(S),\gamma(S))\overset{D_1}{\rightsquigarrow} (M_1,\gamma_1)\overset{D_2}{\rightsquigarrow} (D^2\times [0,1], \partial D^2\times [0,1]).$$
The two indices ``$-$" of $\mc B_{--}$ represent the fact that the boundary orientation of $D_1$ does not match with the orientation of $\ell_1$ and, similarly, the boundary orientation of $D_2$ does not match with the orientation of $\ell_2$.
Let $\mc B'_{--}$ be the branched surface associated with the hierarchy $\mc H_{--}$ before the decomposition along $D_2$. See Figure \ref{fig: pretzel branched surface} for pictures of $\mc B'_{--}$ and $\mc B_{--}$.

\bigskip

\begin{figure}[]
    \centering
    \includegraphics[width=0.6\linewidth]{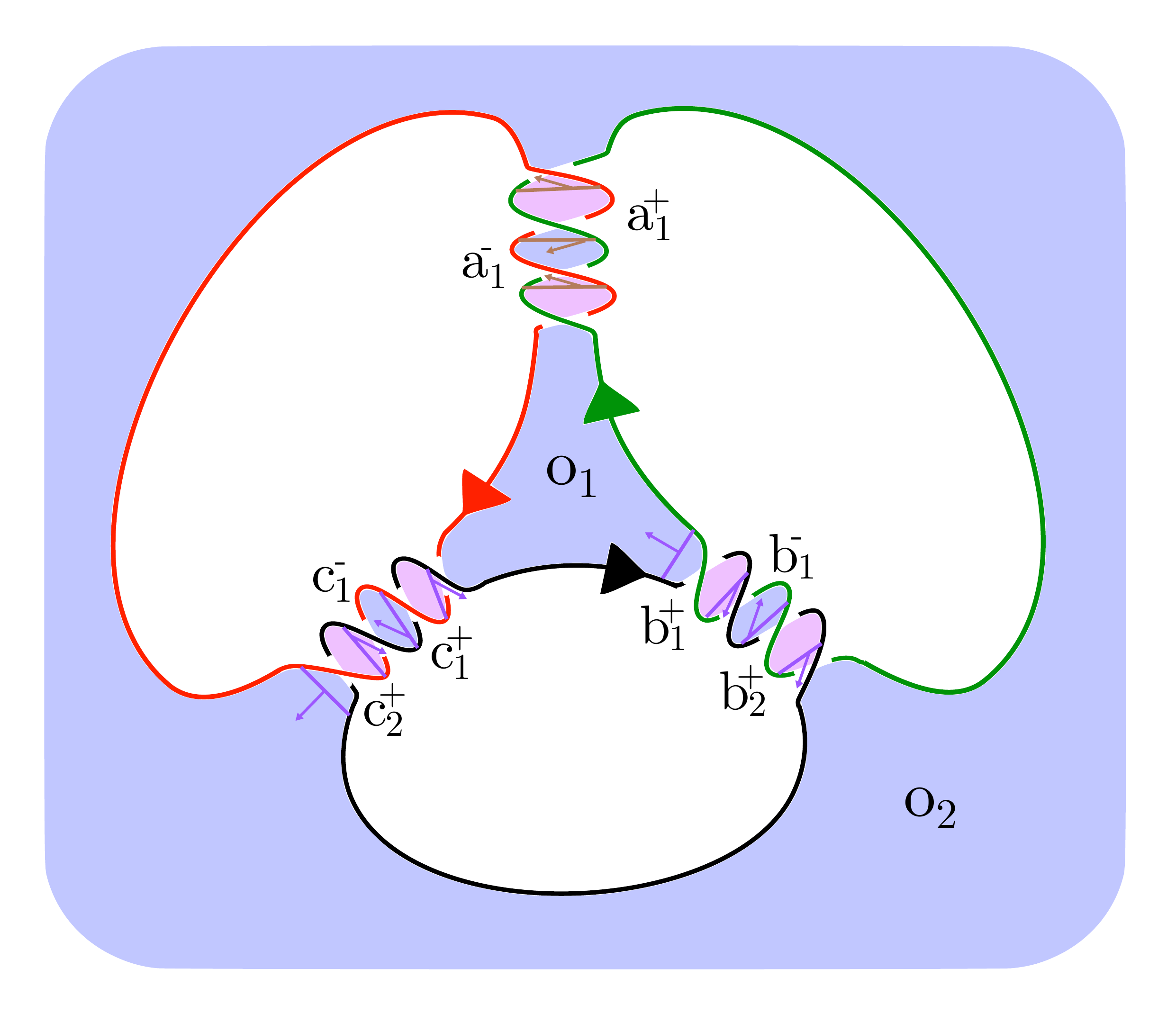}
    \caption{The subdivision of $S$ into sectors of $\mc B_{--}$. The purple and tan segments represent the portion of $\brl(\mc B_{--})$ on $S$. The small arrows indicate the direction of the maw vector field, which is always tangent to $S$. Remember that, by convention, the maw vector field always points outward at the boundary.}
    \label{fig: pretzel sectors}
\end{figure}

\noindent\underline{Analysis of the sectors of $\mc B_{--}$ and computation of the maw dual graph.} Let us analyse the sectors of $\mc B_{--}$. The way $S$ is subdivided into sectors of $\mc B_{--}$ can be inspected through Figure \ref{fig: pretzel sectors}. The two $0$-handles of $S$ give rise to two distinct sectors; call $o_1$ the bounded one and $o_2$ the unbounded one. Both $o_1$ and $o_2$ have six edges alternating on $\partial M$ and on the branched locus $\brl(\mc B_{--})$. The coorientation on the edges of $o_2$ is alternating, hence each corner of it is a double corner and $\chi_m(o_2)=-2$. Instead, one edge on the branched locus of $o_1$ is pointing outward, therefore $\chi_m(o_1)=-1$. The $1$-handle of $S$ corresponding to the $2a$-twisted region is subdivided into $|2a|-2$ square sectors which can be separated in two groups: the edges on the branched locus of the sectors $a_1^+,...,a_{|a|-1}^+$ both point outward, whereas the edges on the branched locus of the sectors $a_1^-,...,a_{|a|-1}^-$ both point inward. 

The sectors $a_i^+$ have $\chi_m=+1$, whilst the sectors $a_j^-$ have $\chi_m=-1$. A similar analysis holds for the $1$-handles associated with the $2b$-twisted and $2c$-twisted region. The only difference is that they are now subdivided into $|2b|-1$ and $|2c|-1$ square sectors respectively, and on each region there is one more sector having $\chi_m=+1$ than the ones having $\chi_m=-1$. Notice that the sum of the $\chi_m$ of the sectors of $S$ has to amount to $\chi(S)=-1$.

Three sectors of $\mc B_{--}$ do not belong to $S$. One is the disk $D_2$, whose maw vector field always points out and therefore $\chi_m(D_2)=+1$. The remaining two sectors patch together to give $D_1$. Indeed, the whole $D_1$ is a sector of the branched surface $\mc B'_{--}$, and its maw vector field always points out. The trace of $\partial D_2$ on $D_1$ splits it into two subdisks $d_0$ and $d_1$. The boundary of one such disk, say $d_1$, still has outward-pointing maw vector field everywhere on the boundary and hence $\chi_m(d_1)=+1$. By direct inspection, or by $\chi_m(d_1)+\chi_m(d_0)=\chi_m(D_1)$, we deduce that $\chi_m(d_0)=0$. 

The maw dual graph $\Gamma_m(\mc B_{--})$ represents a class in $H_1(M_L)\cong (\Z m_1)\oplus (\Z m_2)\oplus(\Z m_3)$, where $m_i$ is the oriented canonical meridian of $\ell_i$ for $i=1,2,3$. Since the exterior of $\mc B_{--}$ consists of a single ball, the graph $\Gamma_m(\mc B_{--})$ has just one vertex, and each sector $s$ of $\mc B_{--}$ yields an oriented dual loop $a(s)$ in $M_L$. For the sectors $s\subset S$, there is always an oriented meridian $m_i$ intersecting $\mc B_{--}$ once, positively and inside $s$; hence $[a(s)]=m_i$. In Table \ref{table}, the total contribution $\chi_m(s)a(s)$ to the maw dual graph is computed for each sector $s$ of $\mc B_{--}$. By summing up the elements in the rightmost column, we obtain \begin{align}\label{equation: B--}[\Gamma_m(\mc B_{--})]=(b+c-1)m_1+(|a|-c-1)m_2+(1-|a|-b)m_3.\end{align}

\begin{table}[]
\centering
\begin{tabular}{@{}lllll@{}}
\toprule
$s$   & how many & $\chi_m(s)$ & $[a(s)]\in H_1(M)$ & total contribution \\ \midrule
$o_1$ & $1$      & $-1$        & $m_2$              &    $-m_2$                \\ 
$o_2$ & $1$      & $-2$        & $m_3$              &    $-2m_3$                \\
$a_i^+$ & $|a|-1$      & $+1$        & $m_2$              &    $+(|a|-1)m_2$                \\
$a_i^-$ & $|a|-1$      & $-1$        & $m_3$              &    $(1-|a|)m_3$                \\ 
$b_i^+$ & $b$      & $+1$        & $m_1$              &    $+bm_1$                \\
$b_i^-$ & $b-1$      & $-1$        & $m_3$              &    $+(1-b)m_3$                \\
$c_i^+$ & $c$      & $+1$        & $m_1$              &    $+cm_1$                \\
$c_i^-$ & $c-1$      & $-1$        & $m_2$              &    $+(1-c)m_2$                \\
$d_1$ & $1$      & $+1$        & $m_2-m_1$              &    $+m_2-m_1$                \\
$D_2$ & $1$      & $+1$        & $m_3-m_2$              &    $+m_3-m_2$                \\\bottomrule
\end{tabular}\caption{Contribution of the sectors of $\mc B_{--}$ to $[\Gamma_m(\mc B_{--})].$}\label{table}
\end{table}

%\bigskip

\noindent\underline{Fully marked surfaces.} Notice that each component of $L$ is unknotted, hence bounds a disk in $S^3$. By looking at the diagram $D$, we can find a disk $D_1$ bounded by $\ell_1$ and intersecting $\ell_2\cup\ell_3$ transversely $b+c$ times. Call $S_1:=D_1\cap M_L$. We can define as well disks $D_2,D_3$ bounded by $\ell_2,\ell_3$ respectively so that the surfaces $S_i:=D_i\cap M_L$ satisfy \begin{align*}
    \chi(S_1)=1-b-c, \\ \chi(S_2)=1-|a|-c, \\ \chi(S_3)=1-|a|-b.
\end{align*}
As the boundary of $S_i$ consists of a longitude of $\ell_i$ and a number of meridians of the other components, we have \begin{align*}
    \langle [\Gamma_m(\mc B_{--})], [-S_1]\rangle=\chi(-S_1) \,\text{ and }\,\langle [\Gamma_m(\mc B_{--})], [S_3]\rangle=\chi(S_3).
\end{align*}
Hence, the classes $\ell_1+\ell_2+\ell_3=[S], -\ell_1=[-S_1]$ and $\ell_3=[S_3]$ belong to a common closed top-dimensional cone of the Thurston norm. The same holds for the classes $-\ell_1-\ell_2-\ell_3, \ell_1$ and $-\ell_3$.

\bigskip

\noindent\underline{Second hierarchy and conclusion of the alternating case.} Notice that we obtain a different groomed hierarchy $\mc H_{-+}$ for $(M_L,\partial M_L)$ if we reverse the orientation of $D_2$. By Lemma \ref{lemma: swapping orientation}, the new branched surface $\mc B_{-+}$ satisfies \begin{align*}[\Gamma_m(\mc B_{-+})]&=[\Gamma_m(\mc B_{--})]-(2-|D_2\cap \gamma_1|)a_{\mc B_{--}}(D_2)\\ &=(b+c-1)m_1+(|a|-c-1)m_2+(1-|a|-b)m_3-(2-2|a|)(m_3-m_2)\\ &=(b+c-1)m_1+(1-|a|-c)m_2+(|a|-b-1)m_3.\end{align*}

We deduce that the classes $\ell_1+\ell_2+\ell_3, -\ell_1$ and $\ell_2$ belong to the same top-dimensional cone of the Thurston norm, and this is true also for $-\ell_1-\ell_2-\ell_3, \ell_1$ and $-\ell_2$

The proof of the alternating case (when $a,b,c$ all have the same sign) now follows by iterating the argument above after a suitable relabelling of the components of $L$. More specifically, one can construct hierarchies analogous to $\mc H_{--}$ and $\mc H_{-+}$ after clockwise rotating the diagram $D$ by $120$ or $240$ degrees. Six hierarchies are hence defined, and the taut foliations associated to them together with their opposites induce the twelve top-dimensional faces of the polyhedron on the left of Figure \ref{fig: polyhedra}. 

\bigskip

\begin{figure}
    \centering
    \includegraphics[width=0.8\linewidth]{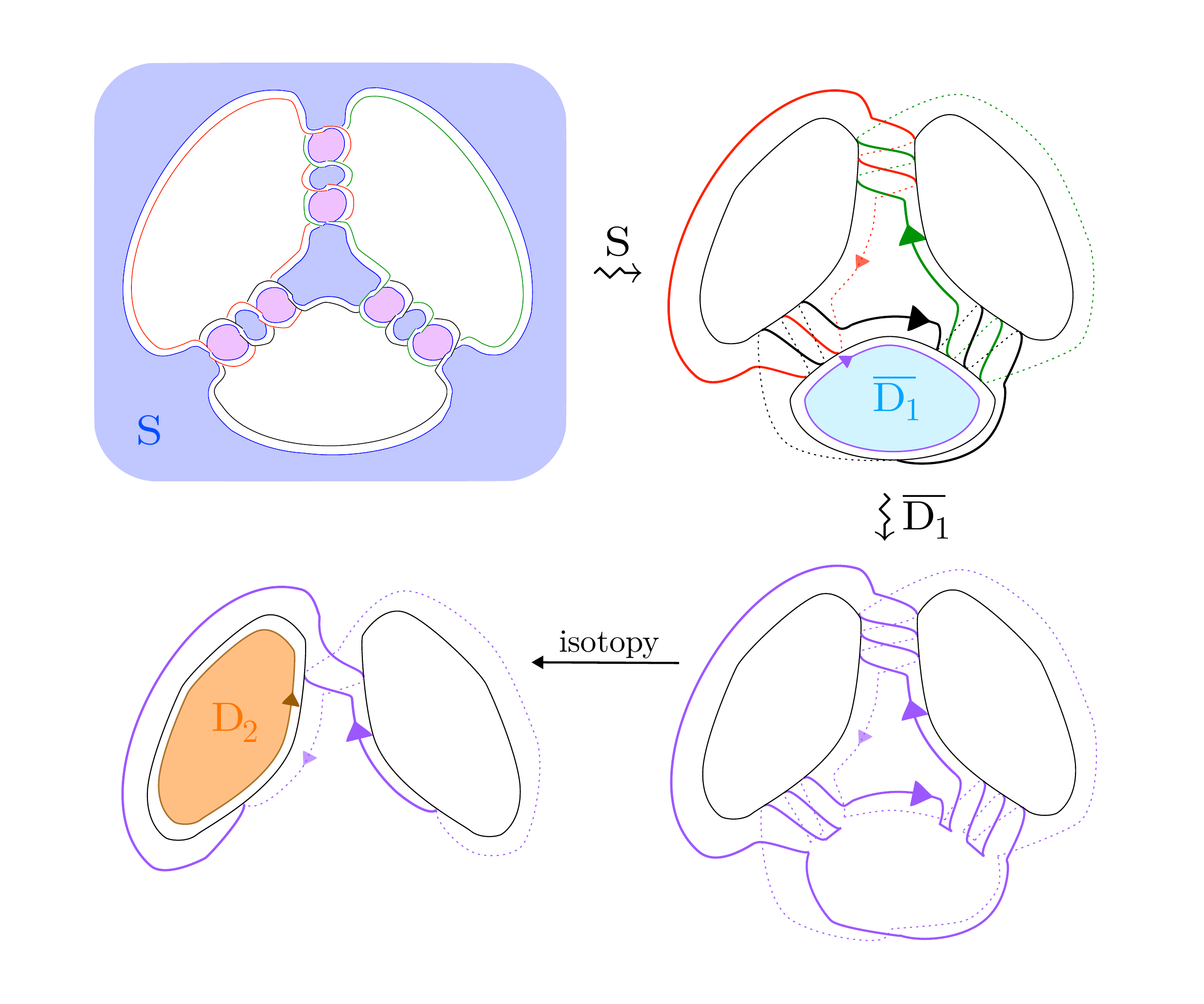}
    \caption{Top left: the Seifert surface $S$ embedded in the exterior of the link $L=P(-4,4,4)$. Top right: the sutured manifold $(M(S),\gamma(S))$ obtained by decomposing $(M_L,\partial M_L)$ along $S$. Bottom right: the sutured manifold $(M_2,\gamma_2)$ resulting by the decomposition of $(M(S),\gamma(S))$ along $\overline{D_1}$. Bottom left: the manifold $(M_2,\gamma_2)$ after an isotopy of $\gamma_2$. After decomposing along $D_2$, the sutured manifold $(M_2,\gamma_2)$ yields a sutured manifold homeomorphic to a ball with one suture, hence a product ball.}
    \label{fig: pretzel hierarchy 2}
\end{figure}

\begin{figure}
    \centering
    \includegraphics[width=0.5\linewidth]{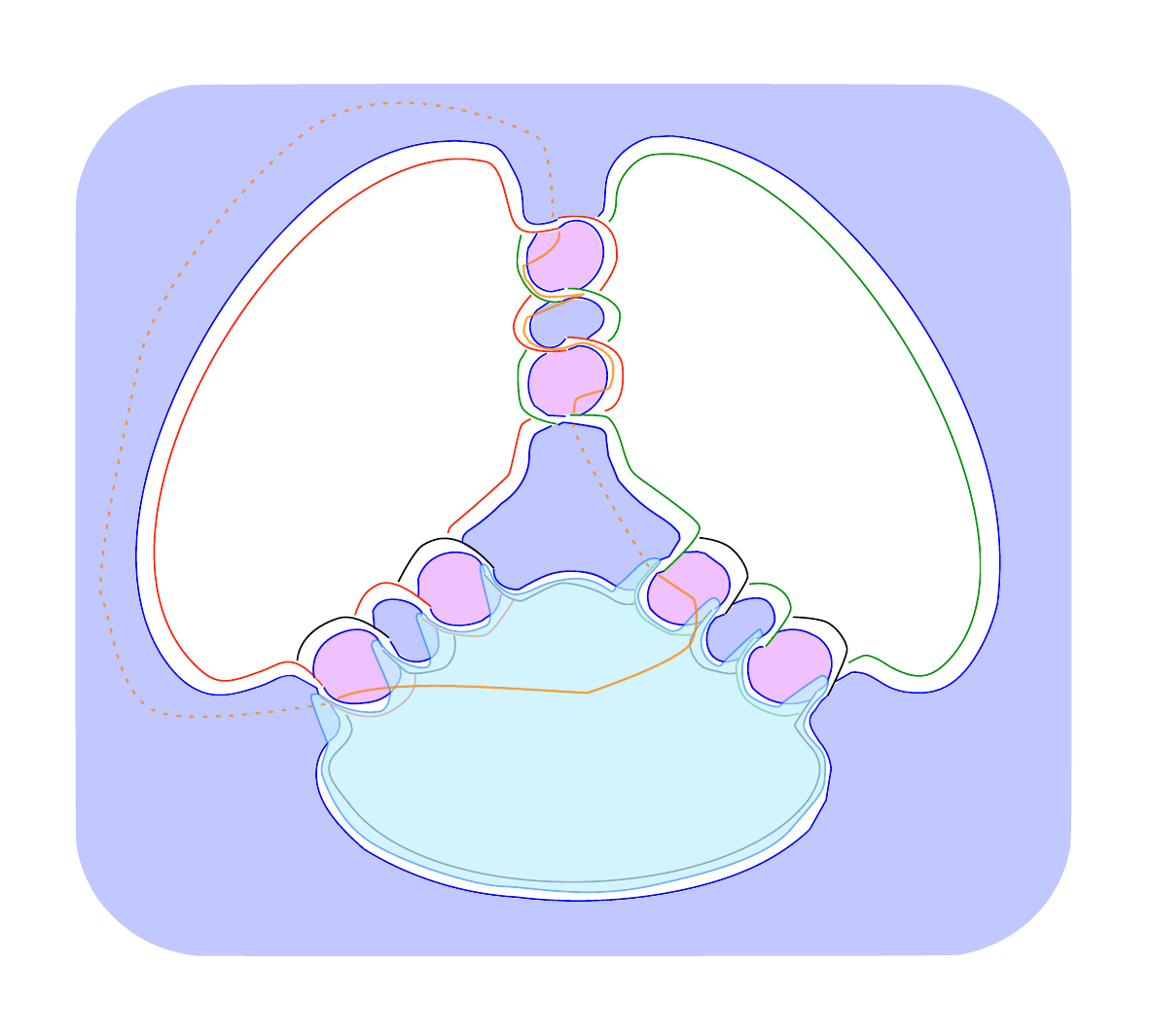}
    \caption{The branched surface $\mc B_{+-}$ for $P(-4,4,4)$. The sector coming from the decomposing disk $D_2$ is not shown for a clearer picture. The trace of $\partial D_2$ is drawn in orange.}
    \label{fig: pretzels branched surface 2}
\end{figure}

\noindent\underline{Case $a\le -2$. Third and fourth hierarchies.} Suppose now that the diagram $D$ is not alternating. If $a\le -2$, consider the groomed sutured manifold hierarchy $\mc H_{+-}$: $$(M_L,\partial M_L)\overset{S}{\rightsquigarrow} (M(S),\gamma(S))\overset{\overline{D_1}}{\rightsquigarrow} (M_2,\gamma_2)\overset{D_2}{\rightsquigarrow} (D^2\times [0,1], \partial D^2\times [0,1]).$$ The case of $P(-4,4,4)$ is depicted in Figure \ref{fig: pretzel hierarchy 2}. Notice that the decomposition along $\overline{D_1}$ reduces by two the number of half-twists of $\gamma_2$ coming from the $2a$-twisted band. The assumption $a\ne -1$ is therefore necessary for $\mc H_{+-}$ to be a taut hierarchy (otherwise $R(\gamma_2)$ would be compressible). Let $\mc B_{+-}$ be the branched surface associated with $\mc H_{+-}$. See Figure \ref{fig: pretzels branched surface 2} for a partial picture of $\mc B_{+-}$ in the exterior of the link $P(-4,4,4)$.

\begin{figure}
    \centering
    \includegraphics[width=0.6\linewidth]{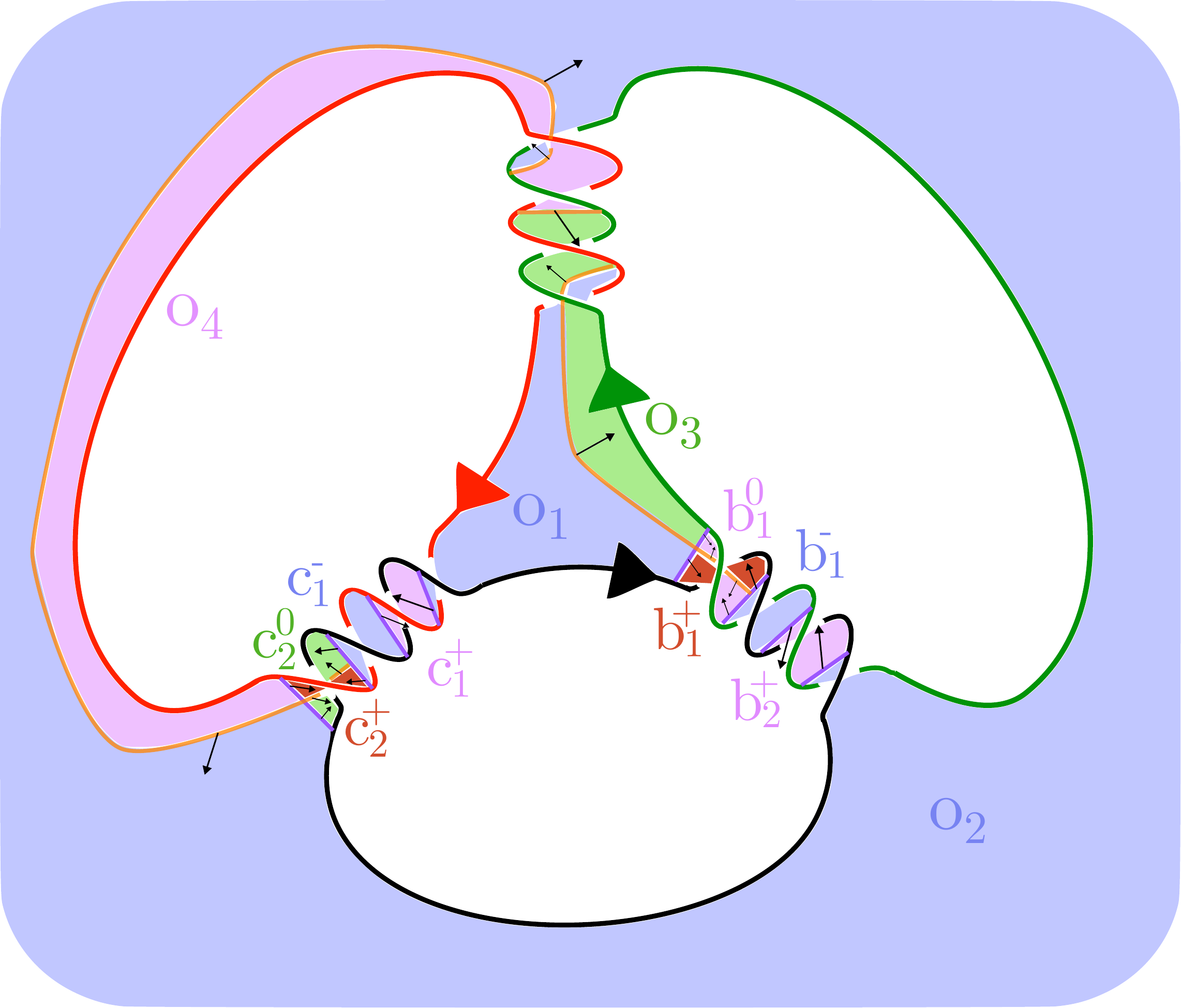}
    \caption{The subdivision of $S$ into sectors of $\mc B_{+-}$. The purple and orange segments represent the portion of $\brl(\mc B_{+-})$ on $S$. The small arrows indicate the direction of the maw vector field, which is always tangent to $S$. Different sides of the same sector are coloured likewise in the figure. When $|a|>2$, $2|a|-4$ more sectors appear on the $2a$-twisted strand. These sectors alternate between sinks ($a_i^+$) and sources ($a_i^-$), see the table below.}
    \label{fig: pretzel sectors two}
    \vspace{10pt}
\end{figure}

A totally analogous analysis of the sectors of $\mc B_{+-}$ as for the previous hierarchy yields the following computation:
$$[\Gamma_m(\mc B_{+-})]=(1-b-c)m_1+(|a|+c-3)m_2+(1+b-|a|)m_3.$$
See Figure \ref{fig: pretzel sectors two} and Table \ref{tablepretzels2} at the end of the paper for further details. If we decompose $(M_2,\gamma_2)$ along $\overline {D_2}$ we get a hierarchy $\mc H_{++}$ whose associated branched surface $B_{++}$ has
\begin{align*}[\Gamma_m(\mc B_{++})]&=[\Gamma_m(\mc B_{+-})]-(2-|D_2\cap \gamma_2|)a_{\mc B_{+-}}(D_2)\\ &=(1-b-c)m_1+(|a|+c-3)m_2+(1+b-|a|)m_3-(4-2|a|)(m_3-m_2)\\ &=(1-b-c)m_1+(1+c-|a|)m_2+(|a|+b-3)m_3.\end{align*}

\bigskip

\noindent\underline{Conclusion of the case $a\le -2$.} Thurston's Inequality guarantees that the Thurston polyhedron $B_x$ of $H_2(M_L,\partial M_L;\R)$ is contained in the polyhedron \begin{align}\label{polyhedron B} B:=\bigcap_{i,j\in \{+,-\}}\{\alpha\in H_2(M_L,\partial M_L;\R)\,|\; |\langle [\Gamma_m(\mc B_{ij})],\alpha\rangle|\le 1\}.\end{align}
By direct inspection -- or by means of a graphic calculator -- the reader can verify that $B$ is the polyhedron in Figure \ref{fig: polyhedra} right, with vertices $$\pm \frac {\ell_1}{b+c-1},\pm \frac {c\ell_1+(b+c-1)\ell_2}{(|a|-1)(b+c-1)},\pm \frac {b\ell_1+(b+c-1)\ell_3}{(|a|-1)(b+c-1)}, \text{ and } \pm (\ell_1+\ell_2+\ell_3).$$

In Lemma \ref{lemma: surface Q}, we show that the classes $c\ell_1+(b+c-1)\ell_2$ and $b\ell_1+(b+c-1)\ell_3$ have Thurston norm $bc+1-b-c+|a-ab-ac-bc|$. If $a\ge ab+bc+ac$, we have $$x(c\ell_1+(b+c-1)\ell_2)=1-b-c+a-ab-ac=(-1-a)(b+c-1)=(|a|-1)(b+c-1).$$
We conclude that in this case the Thurston ball $B_x$ coincides with the polyhedron $B$. In Table \ref{table2} are the Euler classes associated with the eight faces of the Thurston ball $B_x$ for pretzel links $P(2a,2b,2c)$ with $a\ge ab+bc+ac$, $a\le -2$ and $b,c\ge 0$.
\ep

\begin{table}[h]
\centering
\begin{tabular}{c c}
\toprule
Euler class of the face  &  algebraically fully marked classes \\ \midrule
$[\Gamma_m(\mc B_{--})]$ & $\ell_1+\ell_2+\ell_3,-\ell_1, -c\ell_1-(b+c-1)\ell_2, b\ell_1+(b+c-1)\ell_3$ \\ 
$[\Gamma_m(\mc B_{-+})]$ & $\ell_1+\ell_2+\ell_3,-\ell_1, c\ell_1+(b+c-1)\ell_2, -b\ell_1-(b+c-1)\ell_3$  \\ 
$[\Gamma_m(\mc B_{+-})]$ & $\ell_1+\ell_2+\ell_3,\ell_1, b\ell_1+(b+c-1)\ell_3$     \\ 
$[\Gamma_m(\mc B_{++})]$ & $\ell_1+\ell_2+\ell_3,\ell_1, c\ell_1+(b+c-1)\ell_2$    \\ \midrule
$-[\Gamma_m(\mc B_{--})]$ & $-\ell_1-\ell_2-\ell_3,\ell_1, c\ell_1+(b+c-1)\ell_2, -b\ell_1-(b+c-1)\ell_3$ \\ 
$-[\Gamma_m(\mc B_{-+})]$ & $-\ell_1-\ell_2-\ell_3,\ell_1, -c\ell_1-(b+c-1)\ell_2, b\ell_1+(b+c-1)\ell_3$  \\ 
$-[\Gamma_m(\mc B_{+-})]$ & $-\ell_1-\ell_2-\ell_3,-\ell_1, -b\ell_1-(b+c-1)\ell_3$     \\ 
$-[\Gamma_m(\mc B_{++})]$ & $-\ell_1-\ell_2-\ell_3,-\ell_1, -c\ell_1-(b+c-1)\ell_2$     \\\bottomrule
\end{tabular}\caption{Faces of the Thurston ball in the case $a\ge ab+bc+ac$, $a\le -2$ and $b,c\ge 0$.}\label{table2}
\end{table}
%\noindent\underline{Case $a=-1$.} We now outline the proof of the remaining case, leaving the detailed steps to the interested reader. In Lemma \ref{lemma: annulus}, we show that $x(\ell_1+\ell_2)=b-1$ and $x(\ell_1+\ell_3)=c-1$. As $(\ell_1+\ell_2)+(-\ell_1-\ell_3)+(\ell_1+\ell_2+\ell_3)=\ell_1+2\ell_2$, the construction of Lemma \ref{lemma: surface Q} can be reproduced to find a taut surface representing $\ell_1+2\ell_2$ and show that the classes $\ell_1+\ell_2$, $-(\ell_1+\ell_3)$ and $\ell_1+\ell_2+\ell_3$ lie on the same cone. A groomed sutured manifold hierarchy for $(M_L,\partial M_L)$ starting from the annulus $A$ of Lemma \ref{lemma: annulus} and symmetry arguments conclude the proof.

To complete the proof of Theorem \ref{thm: 3-pretzel norm}, it remains to establish Lemma \ref{lemma: surface Q}. This lemma proves the existence of a taut surface in $M_{\ell_1\cup\ell_2}$ that remains norm-minimising after being restricted to $M_L$. This property will be crucial in the proof of Theorem \ref{thm: wrapping number}.

\bl\label{lemma: surface Q} Let $L$ be the pretzel link $P(2a,2b,2c)$ with oriented components $\ell_1,\ell_2$ and $\ell_3$ as in Figure \ref{fig: pretzels}. Suppose $a<0$ and $b,c>0$. Then, there is a taut surface $Q$ in the exterior of $\ell_1\cup \ell_2$ such that:
\begin{itemize}
    \item $Q$ represents the class $c\ell_1+(b+c-1)\ell_2\in H_2(M_{\ell_1\cup\ell_2},\partial M_{\ell_1\cup\ell_2})$, and
    \item $Q$ intersects $\ell_3$ always with the same sign.
\end{itemize}
In particular, the Thurston norm of $c\ell_1+(b+c-1)\ell_2$ in $H_2(M_L,\partial M_L)$ is $$x(c\ell_1+(b+c-1)\ell_2)=bc+1-b-c+|a-ab-ac-bc|.$$
Similarly, $$x(b\ell_1+(b+c-1)\ell_3)=bc+1-b-c+|a-ab-ac-bc|.$$
\el
\bp Notice that the class $\ell_1+\ell_2\in H_2(M_{\ell_1\cup\ell_2},\partial M_{\ell_1\cup\ell_2})$ is represented by an oriented annulus $A$, and the class $\ell_2\in H_2(M_{\ell_1\cup\ell_2},\partial M_{\ell_1\cup\ell_2})$ is represented by a disk $P$ pierced $c$ times by $N(\ell_1)$. In other words, $P$ is a planar surface whose boundary components are a longitude of $\ell_2$ and $c$ meridians of $\ell_1$. Let $Q_0\subset M_{\ell_1\cup\ell_2}$ be the oriented cut and paste of $c$ parallel copies of $A$ with $b-1$ parallel copies of $P$. The surface $Q_0$ represents the class $c\ell_1+(b+c-1)\ell_2\in H_2(M_{\ell_1\cup\ell_2},\partial M_{\ell_1\cup\ell_2})$ and $$\chi(Q_0)=c\chi(A)+(b-1)\chi(P)=(b-1)(1-c)=b+c-1-bc.$$

We claim that there is a surface $Q$ that is isotopic to $Q_0$ in $M_{\ell_1\cup\ell_2}$ and that intersects $\ell_3$ always with the same sign. 
Let us now explain how the second part of the lemma follows from this claim. Since $Q$ is norm-minimising and always intersects $\ell_3$ with the same sign, $Q\cap M_L$ is norm-minimising too. Moreover, the number of intersections between $\ell_3$ and $Q$ is the absolute value of their algebraic intersection number:
\begin{align*}
    \langle [Q],[\ell_3]\rangle&=\langle c[P_1]+(b+c-1)[P],[\ell_3]\rangle\\ &=c\langle [P_1],[\ell_3]\rangle+(b+c-1)\langle[P],[\ell_3]\rangle\\ &=-bc-(b+c-1)a\\ &=a-ab-ac-bc,
\end{align*}
where $P_1\subset M_{\ell_1\cup\ell_2}$ is a punctured disk representing the class $\ell_1$. We conclude that the surface $Q\cap M_L$ is a norm-minimising representative of the class $c\ell_1+(b+c-1)\ell_2\in H_2(M_L,\partial M_L)$ and $$\chi(Q\cap M_L)=\chi(Q)-|Q\cap \ell_3|=b+c-1-bc-|a-ab-ac-bc|,$$ as wanted. The statement about the Thurston norm of $b\ell_1+(b+c-1)\ell_3$ follows from a symmetry between $\ell_2$ and $\ell_3$ in the arguments.

\begin{figure}
    \centering
    \includegraphics[width=0.5\linewidth]{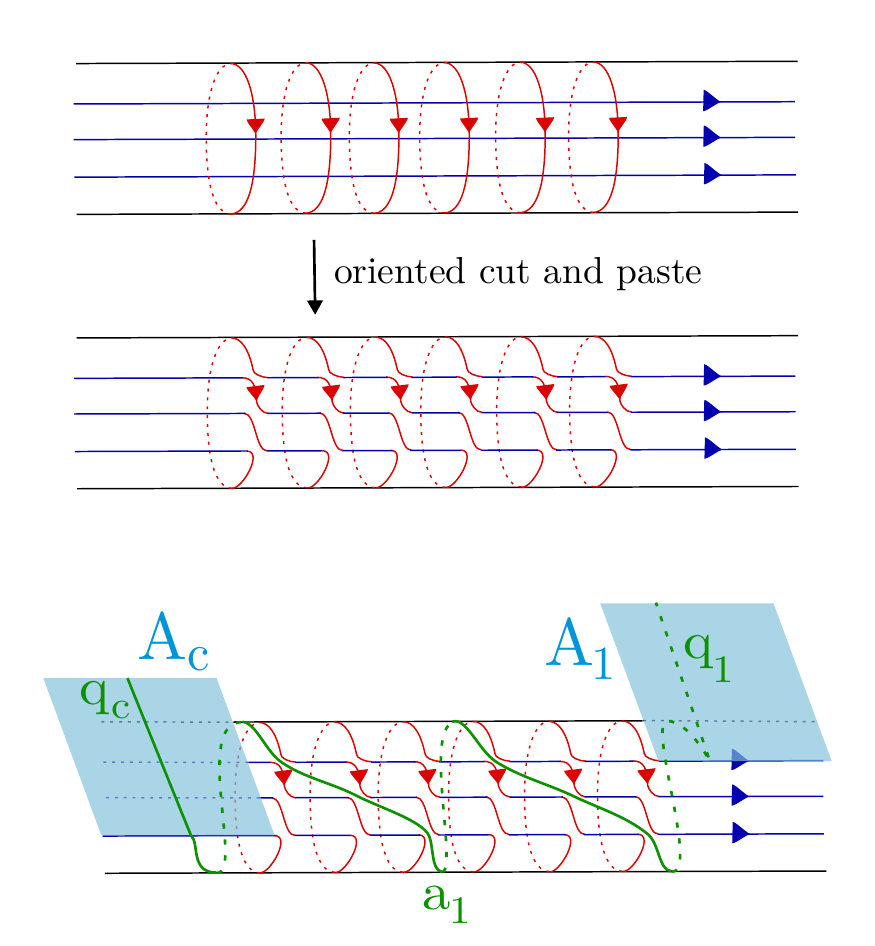}
    \caption{Top: boundary components of $c$ copies of $A$ (in blue) and $b-1$ copies of $P$ (in red) on $\partial N(\ell_1)$. In the picture, $b=c=3$. Middle: boundary components of $Q_0$ on $\partial N(\ell_1)$. Bottom: how to isotope $\ell_3$ (in green) so to intersect $Q_0$ only close to $\partial N(\ell_2)$.}
    \label{fig: boundary intersections}
\end{figure}

\bigskip

We are left to show that $\ell_3$ can be isotoped to intersect $Q_0$ always with the same sign. Remember that $Q_0$ is obtained as the oriented cut and paste of $c$ copies of the annulus $A$ with $b-1$ copies of the planar surface $P$. Let $A_1,...,A_c$ be the parallel copies of $A$ in some circular order coherent with the oriented meridian $m_1$ of $\ell_1$. In other words, if we start walking along $m_1$ from $m_1\cap A_1$, we meet in exact order $A_2$, then $A_3$ and so on. Consider each annulus $A_i$ as interval-fibered over the circle.
The boundary of $P$ on $\partial N(\ell_1)$ consists of $c$ copies of $m_1$. Hence, with respect to the (non-canonical) longitude $\lambda_1':=\partial A\cap N(\ell_1)$ and to the oriented meridian $m_1$ of $\ell_1$, $\partial Q_0$ consists of $c$ parallel simple closed curves representing $\lambda_1'+(b-1)m_1\in H_1(\partial N(\ell_1))$. As a consequence (Figure \ref{fig: boundary intersections}), we can isotope $\ell_3$ in $M_{\ell_1\cup\ell_2}$ so that $\ell_3$ is the union of four closed arcs $a_1,a_2,q_1$ and $q_c$ so that: 
\begin{itemize}
    \item $a_1\subset \partial N(\ell_1)$ and the interior of $a_1$ does not intersect $Q_0$;
    \item $q_1$ is an interval fiber of $A_1$ contained in $Q_0$;
    \item $a_2\subset \partial N(\ell_2)$;
    \item $q_c$ in an interval fiber of $A_c$ contained in $Q_0$.
\end{itemize} Now, $\ell_3$ can be slightly pushed out of a neighbourhood of $Q_0\cup \partial M_{\ell_1\cup\ell_3}$ so to intersect $Q_0$ only close to $\partial N(\ell_2)$. After a further isotopy, $\ell_3$ will therefore always intersect $Q_0$ with the same sign.
\ep

The proof of Theorem \ref{thm: 3-pretzel norm} is now complete. In a similar fashion to the previous lemma, we can establish the following result, which will be useful to compute the wrapping number of the links in Section \ref{sec: wrapping number}.

\bl\label{lemma: annulus} Let $L$ be the pretzel link $P(2a,2b,2c)$ with oriented components $\ell_1,\ell_2$ and $\ell_3$ as in Figure \ref{fig: pretzels}. Suppose $a<0$ and $b,c>0$. Then, there is a Seifert annulus $A$ for the link $\ell_1\cup\ell_2$ such that the surface $A\cap M_L$ is a norm-minimising representative of $\ell_1+\ell_2\in H_2(M_L,\partial M_L)$. In particular, the Thurston norm of $\ell_1+\ell_2$ in $H_2(M_L,\partial M_L)$ is $$x(\ell_1+\ell_2)=|a|+b-2.$$
Similarly, $$x(\ell_1+\ell_3)=|a|+c-2.$$
\el
\bp By looking at the diagram on the left of Figure \ref{fig: pretzels}, one can find an annulus $A_0$ that spans $\ell_1\cup\ell_2$ and intersects $|a|+b$ times $\ell_3$. Since $a$ and $b$ have opposite signs, one can slightly push $A_0$ so that $\ell_3$ intersects the resulting surface twice less. Let $A$ be the annulus $A_0$ after this isotopy. The surface $A\cap M_L$ has Euler characteristic $2-|a|-b$. By identity (\ref{equation: B--}) in the proof of Theorem \ref{thm: 3-pretzel norm}, we obtain: $$\langle [\Gamma(\mc B_{--})], \ell_1+\ell_2\rangle=|a|+b-2,$$ therefore $A\cap M_L$ is norm-minimising. The statement about $\ell_1+\ell_3$ follows by symmetry.
\ep

%\bq What is the Thurston ball of $H_2(M_L,\partial M_L;\R)$, when $L\subset S^3$ is a pretzel link $P(2a,2b,2c)$ with $a\le ab+bc+ac$, $a\le -2$ and $b,c> 0$?
%\eq
%If $a\le ab+bc+ac$, the classes $\pm (c\ell_1+(b+c-1)\ell_2)$ and $\pm( b\ell_1+ (b+c-1)\ell_3)$ are not algebraically fully marked with any of the Euler classes in the left column of Table \ref{table2}. In particular, the unit ball $B_x$ is now strictly contained in the polyhedron $B$ on the right of Figure \ref{fig: polyhedra}. However, it is easy to fit the occurrences of $\pm (\ell_1+\ell_2)$ and $\pm (\ell_1+\ell_3)$ in Table \ref{table2}. Hence, what is left to compute is the Thurston norm in the cone generated by the classes $\ell_1+\ell_2+\ell_3, \ell_1+\ell_2, -\ell_1-\ell_3$ and $-\ell_1-\ell_2-\ell_3$. 

\section{The wrapping number does not satisfy the triangle inequality}\label{sec: wrapping number}

Let $M$ be a compact oriented $3$-manifold and let $K\subset M$ be a knot. Given an integral homology class $\alpha\in H_2(M,\partial M)$, the \emph{wrapping number} of $\alpha$ with respect to $K$ is defined as $$\wrap_K(\alpha):=\min_{[S]=\alpha}|S\cap K|,$$ where the minimum is taken over all norm-minimising incompressible representatives of $\alpha$ transverse to $K$. The function $\wrap_K$ is nonnegative and satisfies $\wrap_K(n\alpha)=|n|\wrap_K(\alpha)$ for every $n \in \Z$ and $\alpha\in H_2(M,\partial M)$. Therefore, it can be extended linearly to a function $\wrap_K:H_2(M,\partial M;\Q)\to \Q$. If $\wrap_K$ satisfied the triangle inequality, then it could be continuously extended to a seminorm on $H_2(M,\partial M;\R)$. We notice that this is the case within a closed cone of the Thurston norm:

\bl Let $M$ be a compact oriented irreducible and $\partial$-irreducible $3$-manifold. If $\alpha,\beta\in H_2(M,\partial M)$ belong to the same closed top-dimensional cone of the Thurston norm, then $$\wrap_K(\alpha+\beta)\le \wrap_K(\alpha)+\wrap_K(\beta)$$ for any knot $K\subset M$.
\el
\bp Let $S$ and $S'$ be norm-minimising incompressible representatives realising the wrapping numbers $\wrap_K(\alpha)$ and $\wrap_K(\beta)$ respectively. As $S$ and $S'$ attain the minimum number of intersections with $K$, we can discard null-homologous components of $S$ and $S'$ until the remaining surfaces are incompressible, $\partial$-incompressible and norm-minimising representatives realising the wrapping number with respect to $K$ (see Remark \ref{rem: taut surfaces}). 

As in the proof of Corollary \ref{cor: norm detection criterion}, $S$ and $S'$ can now be isotoped so that each component of $S\cap S'$ is essential in both $S$ and $S'$. After possibly a further isotopy, $K$ does not pass through a small neighbourhood of $S\cap S'$. Therefore, the oriented cut and paste $\hat S$ of $S$ and $S'$ satisfies \begin{align*}
    \chi_-(\hat S)&=\chi_-(S)+\chi_-(S') \text{ \;and}\\ |\hat S\cap K|&=|S\cap K|+|S'\cap K|.
\end{align*} Since $\alpha$ and $\beta$ belong to the same cone of the Thurston norm, $\hat S$ is norm-minimising. If a component $B$ of $\hat S$ is compressible, then it is null-homologous. As discarding components of $\hat S$ decreases the number of intersections with $K$, we conclude that $$\wrap_K(\alpha+\beta)\le |\hat S\cap K|=|S\cap K|+|S'\cap K|=\wrap_K(\alpha)+\wrap_K(\beta).$$
\ep

In \cite[Question 2.1]{10.4310/jdg/1563242469}, Baker and Taylor asked whether $\wrap_K$ satisfies the triangle inequality in general, suggesting that the answer is likely negative. The previous computations on pretzel links now confirm this suspicion.

\bt\label{thm: wrapping number} Let $L_0\subset S^3$ be a torus link $T(2,2c)$, for some $c>1$, and let $M_0:=M_{L_0}$ be its exterior. There exist infinitely many knots $K\subset M_0$ such that $\wrap_K$ does not satisfy the triangle inequality in $H_2(M_0,\partial M_0)$. 
\et

In order to prove this theorem, we will need a criterion for computing the wrapping number of a homology class. The following lemma serves this purpose.

\bl\label{lemma: wrap criterion} Let $M$ be a compact oriented irreducible and $\partial$-irreducible $3$-manifold. Let $K\subset M$ be a knot, $M(K):=M-N^\circ(K)$ its exterior, and $S$ a norm-minimising incompressible surface in $M$, transverse to $K$, and without sphere or disk components. If $S\cap M(K)$ is norm-minimising in $M(K)$, then $$\wrap_K([S])=|S\cap K|=x_{M(K)}([S\cap M(K)])-x_M([S]),$$ where $x_{M(K)}$ and $x_M$ indicate the Thurston norms on $H_2(M(K),\partial M(K))$ and $H_2(M,\partial M)$ respectively.
\el
\bp Call $\alpha:=[S]\in H_2(M,\partial M)$ and let $T$ be any norm-minimising incompressible representative of $\alpha$. If $T$ intersects $K$ transversely, \begin{equation}\label{intersection1}
    |T\cap K|=\chi(T)-\chi(T\cap M(K)).
\end{equation} If $T$ has no sphere nor disk components, we can rewrite identity (\ref{intersection1}) as \begin{equation}\label{intersection2}
    |T\cap K|=\chi_-(T\cap M(K))-\chi_-(T)=\chi_-(T\cap M(K))-x_M(\alpha).
\end{equation} Now, notice that the homology class of $T\cap M(K)$ is completely determined by $\alpha$. Indeed, $T\cap M(K)$ represents the Poincar\'e dual of the cohomology class $\alpha_K:=\iota^*(\operatorname{PD}(\alpha))\in H^1(M(K))$, where $\iota:M(K)\to M$ is the inclusion and $\operatorname{PD}$ is the Poincar\'e duality. Hence, the identities in (\ref{intersection2}) imply: \begin{equation}\label{intersection3}
    |T\cap K|\ge x_{M(K)}(\alpha_K)- x_M(\alpha),
\end{equation} where equality holds if $T\cap M(K)$ is a norm-minimising representative of $\alpha_K$, and $T$ has no sphere or disk components, e.g., for $T=S$. 

Notice that if $T$ had sphere or disk components instead, these would be zero in homology because $M$ is irreducible and $\partial$-irreducible. Moreover, removing these components from $T$ would not increase the number of intersections with $K$. We conclude that inequality (\ref{intersection3}) holds for every norm-minimising incompressible representative of $\alpha$.  Since the RHS of inequality (\ref{intersection3}) does not depend on the chosen norm-minimising incompressible representative $T$ of $\alpha$, we obtain a lower bound for the wrapping number of $\alpha$: $$\wrap_K(\alpha)\ge x_{M(K)}(\alpha_K)- x_M(\alpha).$$ As this lower bound is realised by $|S\cap K|$, the lemma follows.
\ep

\bp[Proof of Theorem \ref{thm: wrapping number}] Let $L$ be a pretzel link $P(2a,2b,2c)$, where $c$ is given by hypothesis, whilst $a$ and $b$ satisfy $$\begin{cases}
    a\le -2 \\
    b\ge 2 \\
    ab+bc+ac\le a
\end{cases}$$ Notice that, for any $b\ge 2$, there are infinitely many choices of $a$ satisfying these conditions. Label and orient the components of $L$ as in Figure \ref{fig: pretzels}. We will think of $L_0$ as the sublink $\ell_1\cup\ell_2$ -- hence $M_L$ as embedded in $M_0$ -- and we choose $K:=\ell_3$. Consider the classes \begin{align*}
    \alpha&:=(b-1)\ell_1\;\; \text{ and }  \;\;\;\beta :=c\ell_1+(b+c-1)\ell_2
\end{align*}
in $H_2(M_0,\partial M_0)$. We will show that $$\wrap_K(\alpha+\beta)>\wrap_K(\alpha)+\wrap_K(\beta).$$
By abuse of notation, we will also call $\alpha$ and $\beta$ the classes $(b-1)\ell_1$ and $c\ell_1+(b+c-1)\ell_2$ in $H_2(M_L,\partial M_L)$.
We now verify that each of the classes $\alpha,\beta$ and $\alpha+\beta$ have a surface representative satisfying the hypotheses of Lemma \ref{lemma: wrap criterion}, for $M=M_0$ and $M(K)=M_L$:
\begin{itemize}
    \item We had defined a norm-minimising representative $S_1$ of $\ell_1\in H_2(M_L,\partial M_L)$ by intersecting a disk $D_1$ spanning $\ell_1$ with $M_L$. Therefore the lemma holds for $\alpha$ by taking $b-1$ parallel copies of $D_1\cap M_0$.
    \item Lemma \ref{lemma: surface Q} guarantees that there is a taut surface $Q$ representing $\beta\in H_2(M_0,\partial M_0)$ such that $Q\cap M_L$ is norm-minimising. Therefore, $\wrap_K(\beta)=x_L(\beta)-x_0(\beta)$.
    \item Finally,  the class $\alpha+\beta\in H_2(M_0,\partial M_0)$ is represented by $b+c-1$ parallel copies of the properly embedded annulus $A$ of Lemma \ref{lemma: annulus}. It follows by the lemma that the union of $b+c-1$ parallel copies of $A\cap M_L$ is norm-minimising in $M_L$. 
\end{itemize}
In particular, we can compute the wrapping numbers of $\alpha,\beta$ and $\alpha+\beta$ in terms of the Thurston norms $x_L$ and $x_0$ on $H_2(M_L,\partial M_L)$ and $H_2(M_0,\partial M_0)$ respectively. We obtain: \begin{align*}
    \wrap_K(\alpha+\beta)-\wrap_K(\alpha)&-\wrap_K(\beta)=\\ &=x_L(\alpha+\beta)-x_0(\alpha+\beta)-x_L(\alpha)+x_0(\alpha)-x_L(\beta)+x_0(\beta)\\ &=[x_L(\alpha+\beta)-x_L(\alpha)-x_L(\beta)]+[x_0(\alpha)+x_0(\beta)-x_0(\alpha+\beta)]\\ &= x_0(\alpha)+x_0(\beta)-x_0(\alpha+\beta)\\ &=x_0(\alpha)+x_0(\beta)\ge x_0(\alpha)=c-1>0,
\end{align*}
where we used that the classes $\alpha$ and $\beta$ belong to the same closed cone of $x_L$ and that $x_0(\alpha+\beta)=0$ as $L_0$ spans an annulus. Alternatively, the quantity $x_0(\alpha)+x_0(\beta)-x_0(\alpha+\beta)$ is positive because the classes $\alpha$ and $\beta$ do not belong to the same top-dimensional cone of $x_0$.
\ep

The proof of Theorem \ref{thm: wrapping number} extends to the following more general criterion to detect cases where the wrapping number fails to satisfy the triangle inequality.

\bprop Let $M$ be a compact oriented irreducible and $\partial$-irreducible $3$-manifold. Let $K\subset M$ be a knot and $M(K)$ its exterior. Let $\alpha$ and $\beta\in H_2(M,\partial M)$ be classes such that:
\begin{itemize}
    \item[(i)] $\alpha$ and $\beta$  do \underline{not} belong to the same closed top-dimensional cone of the Thurston norm of $M$.
    \item[(ii)] The classes $\alpha_K$ and $\beta_K$ obtained by restricting $\alpha$ and $\beta$ to $H_2(M(K),\partial M(K))$ belong to the same top-dimensional closed cone of the Thurston norm of $M(K)$.
    \item[(iii)] The classes $\alpha,\beta$ and $\alpha+\beta$ have norm-minimising incompressible representatives that restrict to norm-minimising representatives of $\alpha_K,\beta_K$  and $\alpha_K+\beta_K$ respectively.
\end{itemize}
Then $$\wrap_K(\alpha+\beta)>\wrap_K(\alpha)+\wrap_K(\beta).$$
\eprop

\begin{figure}
    \centering
    \includegraphics[width=0.4\linewidth]{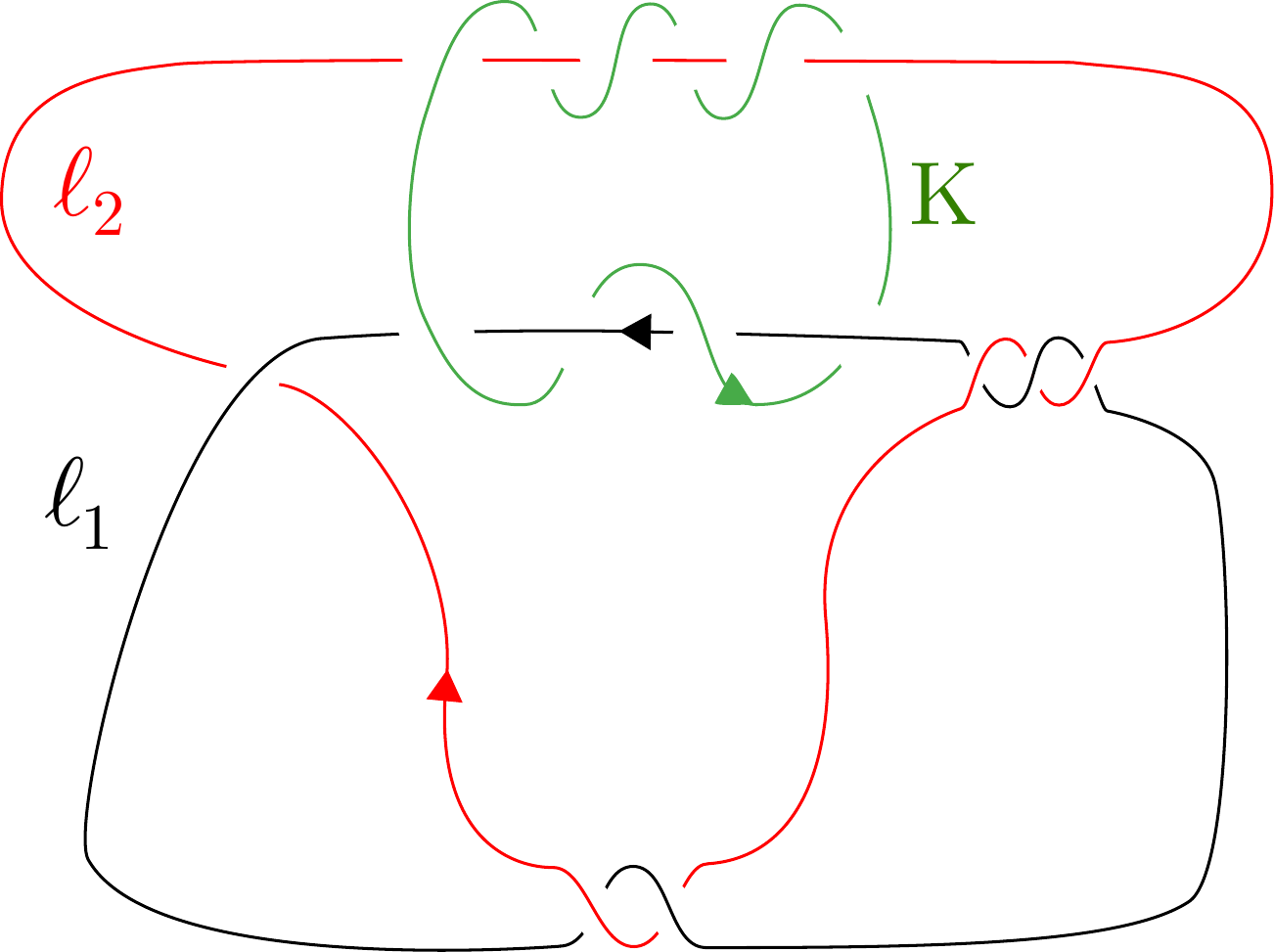}
    \caption{The hyperbolic $2$-bridge link $L_0=\ell_1\cup\ell_2$, and the knot $K\subset M_{L_0}$.}
    \label{fig: hyperbolic wrap}
\end{figure}

\br It is also possible to find hyperbolic $2$-bridge links $L_0$ that satisfy the conclusions of Theorem \ref{thm: wrapping number}. An example is shown in Figure \ref{fig: hyperbolic wrap}. In this case, Lemma \ref{lemma: wrap criterion} applies to the classes $\alpha=\ell_1, \beta=3\ell_1+4\ell_2$ and $\alpha+\beta=4(\ell_1+\ell_2)$. We can then compute that:
\begin{itemize}
    \item $\wrap_K(\alpha)=2$ because a Seifert surface for $\ell_1$ intersects $K$ twice with the same sign;
    \item $\wrap_K(\alpha+\beta)=12$ because a minimal-genus Seifert surface $S_0$ for $L_0$ intersects $K$ three times and remains norm-minimising in $M_{L_0\cup K}$. This statement can be verified by means of a taut sutured manifold hierarchy for $(M_{L_0\cup K},\partial M_{L_0\cup K})$.
    \item $\wrap_K(\beta)=|\langle \beta,K\rangle|=|3\operatorname{Lk}(\ell_1,\ell_3)+4\operatorname{Lk}(\ell_2,\ell_3)|=|-6+12|=6$ thanks to a construction similar to the one in the proof of Lemma \ref{lemma: surface Q}. Indeed, close to $K$, the Seifert surface $S_0$ looks like a rectangle $R=[0,1]\times (0,1)$, with $\{0,1\}\times(0,1)\subset \partial M_0$, and $K\subset N(R\cup \partial M_0)$.
    The construction of Lemma \ref{lemma: surface Q} can be repeated to obtain information about the class $$|\operatorname{Lk} (\ell_1,\ell_2)| \ell_1+(|\operatorname{Lk}(\ell_1,\ell_2)|+|\operatorname{Lk}(\ell_1,\ell_3)|-1)\ell_2=3\ell_1+4\ell_2=\beta.$$
    By \cite{2-bridge_links}, a norm-minimising representative for $\beta$ can be obtained as the oriented cut and paste of three copies of $S_0$ with a norm-minimising representative of $\ell_2$.
\end{itemize}
Also in this case, $\wrap_K(\alpha+\beta)>\wrap_K(\alpha)+\wrap_K(\beta).$
\er

\bibliographystyle{alpha}
\bibliography{bibliography}

\newpage

\begin{table}[]
\centering
\begin{tabular}{@{}lllll@{}}
\toprule
$s$   & how many & $\chi_m(s)$ & $[a(s)]\in H_1(M)$ & total contribution \\ \midrule
$o_1$ & $1$      & $+1$        & $m_2$              &    $+m_2$                \\ 
$o_2$ & $1$      & $0$        & $m_3$             &    $0$                \\
$o_3$ & $1$      & $-1$        & $m_3$              &    $-m_3$                \\
$o_4$ & $1$      & $+1$        & $m_2$              &    $+m_2$                \\
$a_i^+$ & $|a|-2$      & $+1$        & $m_2$              &    $+(|a|-2)m_2$                \\
$a_i^-$ & $|a|-2$      & $-1$        & $m_3$              &    $+(2-|a|)m_3$                \\ 
$b_1^0$ & $1$      & $0$        & $m_1-m_2+m_3$              &    $0$                    \\
$b_i^+$ & $b$      & $-1$        & $m_1$              &    $-bm_1$                \\
$b_i^-$ & $b-1$      & $+1$        & $m_3$              &    $+(b-1)m_3$                \\
$c_i^+$ ($1\le i\le c-1$) & $c-1$      & $-1$        & $m_1$              &    $+(1-c)m_1$                \\
$c_c^+$ & $1$      & $-1$        & $m_1+m_2-m_3$              &    $-m_1-m_2+m_3$                \\
$c_c^0$ & $1$      & $0$        & $m_1$              &    $0$                \\
$c_i^-$ & $c-1$      & $+1$        & $m_2$              &    $+(c-1)m_2$                \\
$d_0$ & $1$      & $+1$        & $m_1-m_3$              &    $+m_1-m_3$                \\
$D_2$ & $1$      & $+1$        & $m_3-m_2$              &    $+m_3-m_2$                \\\bottomrule
\end{tabular}\caption{Contribution of the sectors of $\mc B_{+-}$ to $[\Gamma_m(\mc B_{+-})].$ In order to compute $[a(c_c^+)]$, it might be useful to notice that $[a(c_c^0)]=[a(c_c^+)]+[a(D_2)]$ because of the branching information.}\label{tablepretzels2}
\end{table}

\end{document}